\documentclass{article}
\usepackage{RR}
\usepackage{hyperref}
\usepackage{graphicx,amssymb,amsmath,color}
\usepackage{subfigure}
\usepackage{times}

\newcommand\Rey{\mbox{\textit{Re}}}  

\newcommand{\vct}[1]{\boldsymbol{#1}}
\newcommand{\tond}[1]{\left( #1 \right)}
\RRdate{F\'evrier 2007}
\RRauthor{Marcelo Buffoni\thanks{DIASP - Politecnico di Torino, 10129
    Torino, Italy}\thanks[sfr]{MAB - Universit\'{e} Bordeaux 1 et MC2
    - INRIA Futurs, 33405 Talence, France} 
  \and
  Simone Camarri\thanks[sit]{DIA - Universit\`{a} di Pisa, 56127
    Pisa,Italy}
  \and
  Angelo Iollo\thanksref{sfr}
  \and
  Edoardo Lombardi\thanksref{sfr}\thanksref{sit}
  \and
  Maria Vittoria Salvetti\thanksref{sit}
}
\authorhead{Buffoni \& Camarri \& Iollo \& Lombardi \& Salvetti}
\RRtitle{Un observateur non linéaire pour des écoulements
  tridimensionnels instables} 
\RRetitle{A non-linear observer for unsteady three-dimensional flows}
\titlehead{A non-linear observer for unsteady three-dimensional flows}
\RRresume{Nous proposons une méthode pour estimer le champ de vitesse
  d'un écoulement instable en utilisant un nombre limité de
  mesures. La méthode est basée sur un modèle d'ordre réduit non
  linéaire de l'écoulement et sur l'expansion du champ de vitesse en
  termes de fonctions de base  empiriques. L'idée principale est
  d'imposer que les coefficients de l'expansion modale du champ de
  vitesse donnent la meilleure approximation des mesures disponibles
  et qu'en même temps ils satisfassent aussi étroitement que possible
  le modèle d'ordre réduit non linéaire. L'utilisation pratique peut
  s'étendre du contrôle actif en boucle fermée de l'écoulement à la
  surveillance de l'écoulement dans des régions inaccessibles. La
  technique proposée est appliquée à l'écoulement autour d'un cylindre
  carré confiné, dans des régimes laminaires bi et
  tridimensionnels. Des comparaisons avec des techniques linéaires et
  non linéaires d'évaluation existantes sont fournies.}  
\RRabstract{A method is proposed to estimate the velocity field of an
  unsteady flow using a limited number of flow measurements. The
  method is based on a non-linear low-dimensional model of the flow
  and on expanding the velocity field in terms of empirical basis
  functions. The main idea is to impose that the coefficients of the
  modal expansion of the velocity field give the best approximation to
  the available measurements and that  at the same time they satisfy
  as close as possible the non-linear low-order model. The practical
  use may range from feedback flow control to monitoring of the flow
  in non-accessible regions. The proposed technique is applied to the
  flow around a confined square cylinder, both in two- and
  three-dimensional laminar flow regimes. Comparisons are provided
  with existing linear and non-linear estimation techniques.} 
\RRmotcle{modèles réduits, estimation dynamique, observateurs non
  linéaires}  
\RRkeyword{low-order models, dynamic estimation, non-linear observers}  
\RRprojet{MC2}
\RRtheme{\THNum} 
\URFuturs 
\begin{document}
\makeRR   

\section{Introduction}

The problem of deriving an accurate estimation of the velocity field
in an unsteady complex flow, starting from a limited number of
measurements, is of great importance in many engineering
applications. For instance, in the design of a feedback control, a
knowledge of the velocity field is a fundamental element in deciding
the appropriate actuator reaction to different flow conditions. In
other applications it may be necessary or advisable to monitor the
flow conditions in regions of space which are difficult to access or
where probes cannot be fitted without causing interference problems.

The method that we propose exploits an idea which is similar to that
at the basis of the Kalman filter (see \cite{Kalman}). The starting
point is a Galerkin representation of the velocity field
$\vct{u}(\vct{x},t)$ in terms of $N_r$ empirical eigenfunctions,
$\vct{\Phi}^i(\vct{x})$, obtained by Proper Orthogonal Decomposition
(POD)  (see \cite{L67})
\begin{equation}
\label{galrap}
\vct{u}(\vct{x},t)=\overline{\vct{u}}(\vct{x})+\sum_{i=1}^{N_r}{a_i(t)\vct{\Phi}^i(\vct{x})} 
\end{equation}
where $\vct{u}(\vct{x},t): \mathbb{R}^n \times [0,T] \rightarrow
\mathbb{R}^n$, $\Phi^i(\vct{x}): \mathbb{R}^n \rightarrow
\mathbb{R}^n$, $n \in \{2,3\}$ according to the physical space
dimension, $\overline{\vct{u}}(\vct{x})$ is some reference velocity
field and $a_i(t): I=[0,T] \subset \mathbb{R}  \rightarrow \mathbb{R}$. 

For a given flow, the POD modes can be computed once for all based on
Direct Numerical Simulation (DNS) or on highly resolved experimental
velocity fields, such as those obtained by particle image
velocimetry. An instantaneous velocity field can thus be reconstructed
by estimating the coefficients $a_i(t)$ of its Galerkin representation.

One simple approach to estimating the POD coefficients is to
approximate the flow measurements in a least square sense, as done,
for instance, in \cite{Galletti2004}. 

A similar procedure is also used in the estimation based on gappy POD,
see \cite{karnia} and \cite{willcox}. Another possible approach, the
linear stochastic estimation (LSE), is based on the assumption that a
linear correlation exists  between the flow measurements and the value
of the POD modal coefficients (see, for instance, \cite{Bonnet}). 

However, these approaches encounter difficulties in giving accurate
estimations when three-dimensional flows with complicated unsteady
patterns are considered, or when a very limited number of sensors is
available. Under these conditions, for instance, the least squares
approach cited above (LSQ) rapidly becomes ill conditioned. This
simply reflects the fact that more and more different flow
configurations correspond to the same set of measurements. To
circumvent those problems, many contributions in the literature have
aimed to determine the effective placement of the sensors (see
e.g. \cite{glauser}, \cite{Cohen2004}, \cite{Cohen2006},
\cite{willcox}). For example in \cite{willcox}, a systematic approach
to sensor placement is formulated within the gappy POD framework using
a condition number criterion.

In order to improve estimation performance, extensions of the above
methods have been proposed: quadratic stochastic estimation (QSE)
\cite{adr77}, \cite{naguib2001} and spectral linear stochastic
estimation (SLSE) \cite{ewing99}. They allow more accurate estimations
compared with LSQ or LSE methods, but, in fact, neither of these
methods takes into account the underlying dynamic model that the POD
coefficients must satisfy, i.e., a finite dimensional equivalent of
the Navier-Stokes equations that is obtained by the Galerkin
projection of the flow equations on the POD modes retained for the
representation of the velocity field. In this sense, the aim of the
present study is to discuss an approach that combines a linear
estimation of the coefficients $a_i(t)$ with an appropriate non-linear
low-dimensional flow model. Our objective is not, however, to propose
an estimation method that can be readily implemented for real time
applications, even if a few indications in this direction are
given. Rather, our objective is to understand whether a non-linear
observer outperforms existing linear flow observers, without the
constraints imposed by an actual recursive algorithm, e.g., a
real-time computation. Moreover, instead of what was done, for
example, in \cite{Bewley}, this study is confined to a deterministic
framework, since the model as well as the measurements are supposedly
not affected by noise. If, within this framework, a dynamic estimation
turns out to be less satisfactory than static approaches, i.e., those
which use no model, then there would be little interest in pursuing
the research in this direction. In addition, we address the issue of
the sensitivity of the proposed approach to sensor type and
location. Finally, we present an application to a flow, which is
characterized by a significant three-dimensionality and a non-periodic
dynamics. 

\section{Flow set up and low order model}
\label{sec:flow}

\begin{figure}
  \centering
    \includegraphics[width=10.3cm]{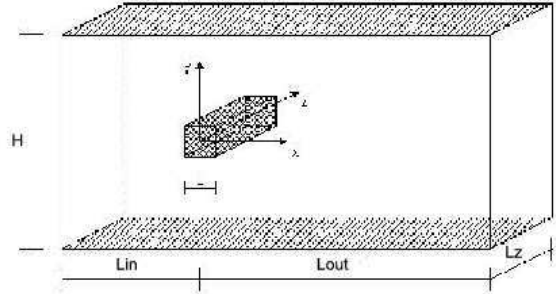}
  \caption{Computational domain $\Omega$.}
  \label{fig:compdom} 
\end{figure}

The flow over an infinitely long square cylinder symmetrically
confined by two parallel planes is considered. A sketch showing the
geometry, the frame of reference and the adopted notation is plotted
in figure~\ref{fig:compdom}. At the inlet, the incoming flow is
assumed to have a Poiseuille profile with maximum center-line velocity
$U_c$.  
Two Reynolds numbers $Re=U_c L / \nu$ were considered, one at which
the flow is two-dimensional ($Re=150$) and the other one leading to a
three-dimensional flow in the wake ($Re=300$). With reference to
figure~\ref{fig:compdom}, $L/H=1/8$, $Lin/L = 12$, $Lout/L = 20$ and $
Lz/L= 0.6$ for the two-dimensional case, whereas $L/Lz = 6$ for the
three-dimensional one. Periodic boundary conditions are imposed in the
span-wise direction and no-slip conditions are enforced both on the
cylinder and on the parallel walls. Details concerning the grids and
the numerical set up are reported in \cite{Buffoni2006}. All the
quantities reported in the following have been made non-dimensional by
$L$ and $U_c$. The two-dimensional flow obtained at $Re=150$ is a
classic vortex street with a well defined shedding frequency. However,
the interaction with the confining walls adds to the complexity of the
flow and leads to some peculiar features, like the fact that the
vertical position of the span-wise vortices is opposite to the one in
the classic von K\'arm\'an street (\cite{Camarri2006}). For the
three-dimensional case the situation is even more complex, due to
instabilities developing in the span-wise direction. The flow is no
longer periodic and exhibits complicated flow patterns
(\cite{Buffoni2006}).  

The POD modes $\vct{\Phi}^k(\vct{x})$ are found  using the snapshot
method (\cite{S87})
\[
\begin{array}{l}
\vct{\Phi}^k=\sum_{i=1}^N{b_i^k \, \vct{U}^{(i)}}
\end{array}
\]
where $\vct{U}^{(i)}=\vct{u}(\vct{x},t_i)$ are flow snapshots taken at
time $t_i \in [0,T]$, $N$ is the number of snapshots, $k \in
\{1,\dots,N\}$, and the coefficients $b_i^k \in \mathbb{R}$ are such
that the vectors $(b_1^k,\dots,b_N^k)$ are the eigenvectors of the
time correlation matrix $\int_\Omega{\vct{U}^{(j)} \cdot \vct{U}^{(l)}
  \, dx}$, of size $N \times N$. Only a limited number of modes,
$N_r$, is used to represent the velocity field. In particular we took
$N_r$ = 6 and $N_r$ = 20 for the two-dimensional and the
three-dimensional cases, respectively. 

A Galerkin projection of the incompressible Navier-Stokes equations
over the retained POD modes has been carried out. This leads to the
following $N_r$-dimensional dynamical system 
\begin{equation}
\label{model}
\begin{array}{l}
R_r(\vct{a}(t))=\dot{a}_r(t) - A_r - C_{kr} a_k (t) + B_{ksr}a_k(t)a_s(t)=0\\ 
a_r (0) = (\vct{u}(\vct{x},0)-\overline{\vct{u}}(\vct{x}),\vct{\Phi}^r )
\end{array}
\end{equation}
where $\vct{a}(t): I \rightarrow \mathbb{R}^{N_r}$ and
$\vct{a}(t)=\{a_1(t),\dots,a_{N_r}(t)\}$; $r$, $k$ and $s$ run from
$1$ to $N_r$ and the Einstein summation convention is used.  The
scalar coefficients $B_{ksr}$ come directly from the Galerkin
projection of the non-linear terms in the Navier-Stokes equations, and
they can easily be expressed in terms of the POD modes. The scalar
terms  $A_r$ and $C_{kr}$ are calibrated using a pseudo-spectral
method  to take into account the pressure drop, as well as the
interaction of unresolved modes in the POD expansion. The calibration
consists in solving an inverse problem, where the coefficients $A_r$
and $C_{kr}$ are found in order to minimize the difference, measured
in the $L^2$ norm, between the model prediction and the actual
reference solution. See \cite{Galletti2006} for a detailed discussion
of the calibration technique. 

\begin{figure}[h!]
  \centering
  \subfigure{
    \includegraphics[width=3cm,height=3cm]{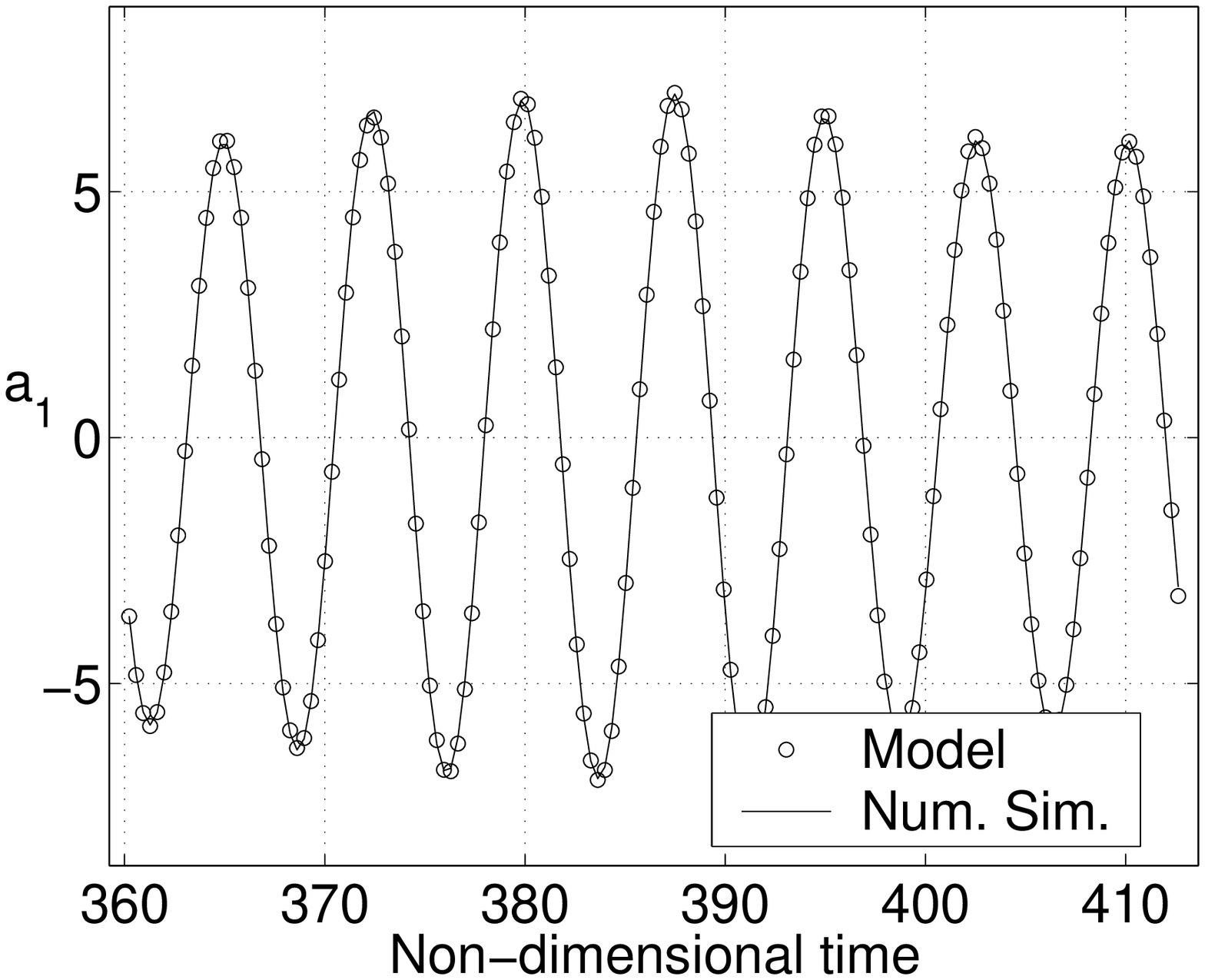}}
  \subfigure{
    \includegraphics[width=3cm,height=3cm]{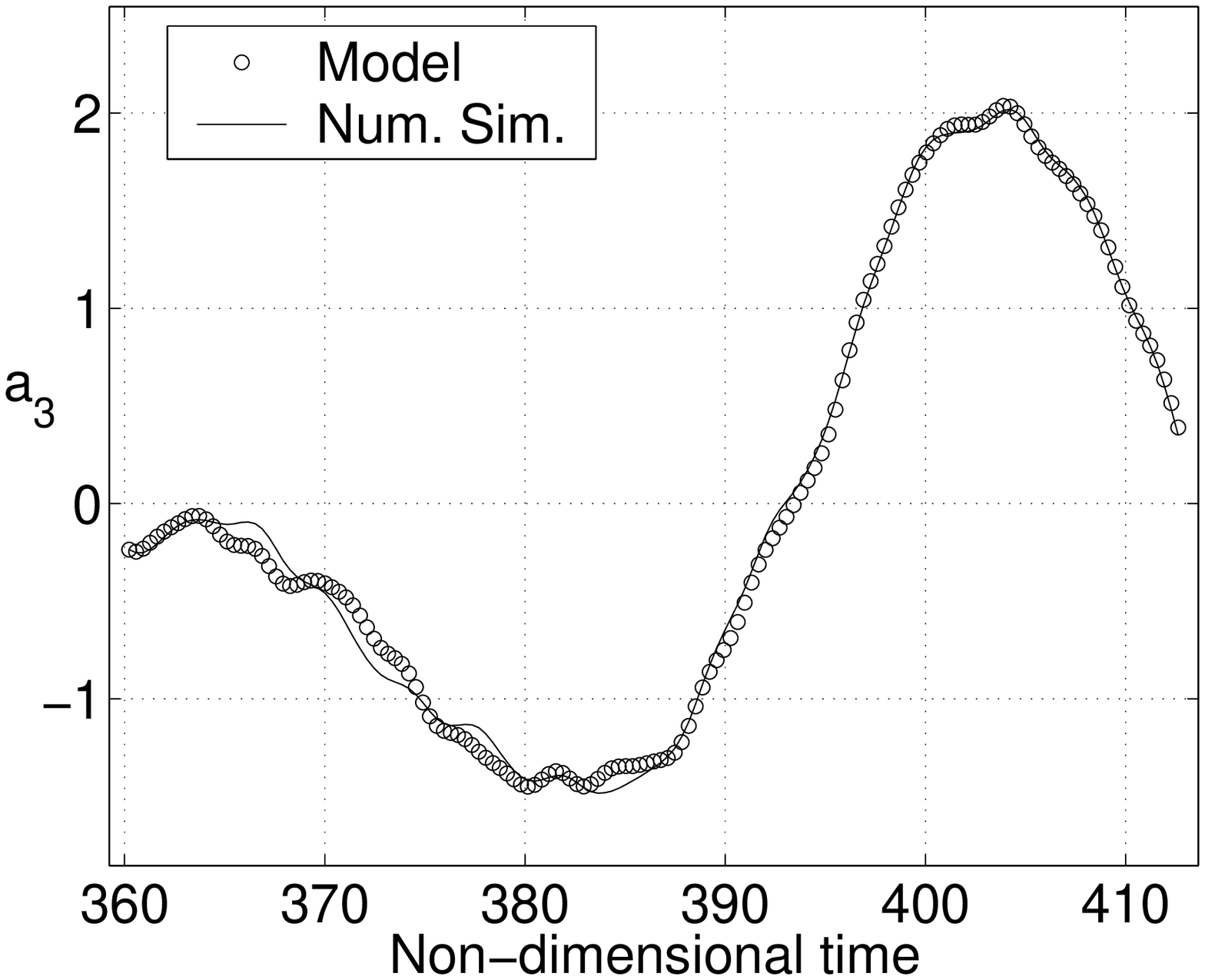}}
  \subfigure{
    \includegraphics[width=3cm,height=3cm]{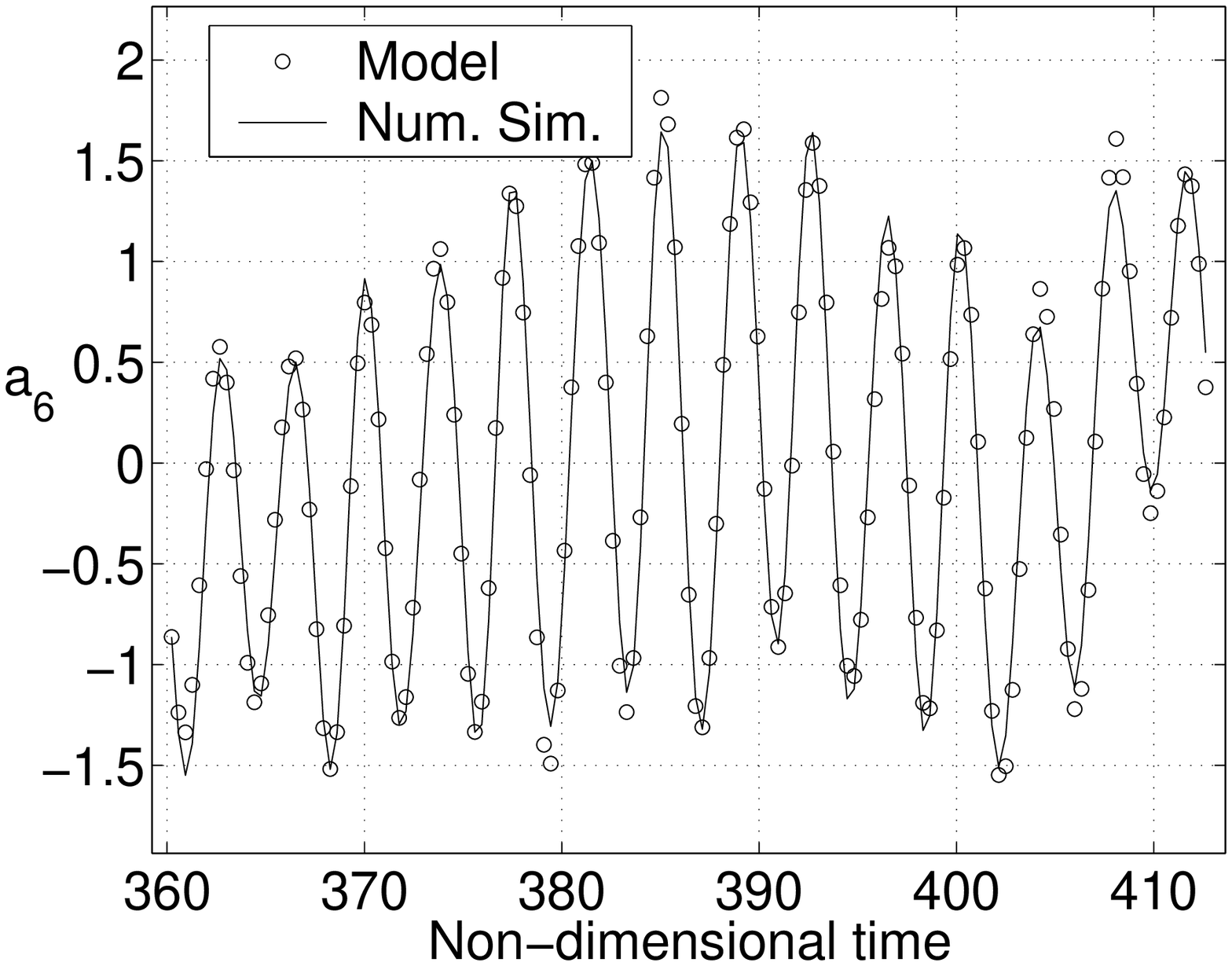}}
  \subfigure{
    \includegraphics[width=3cm,height=3cm]{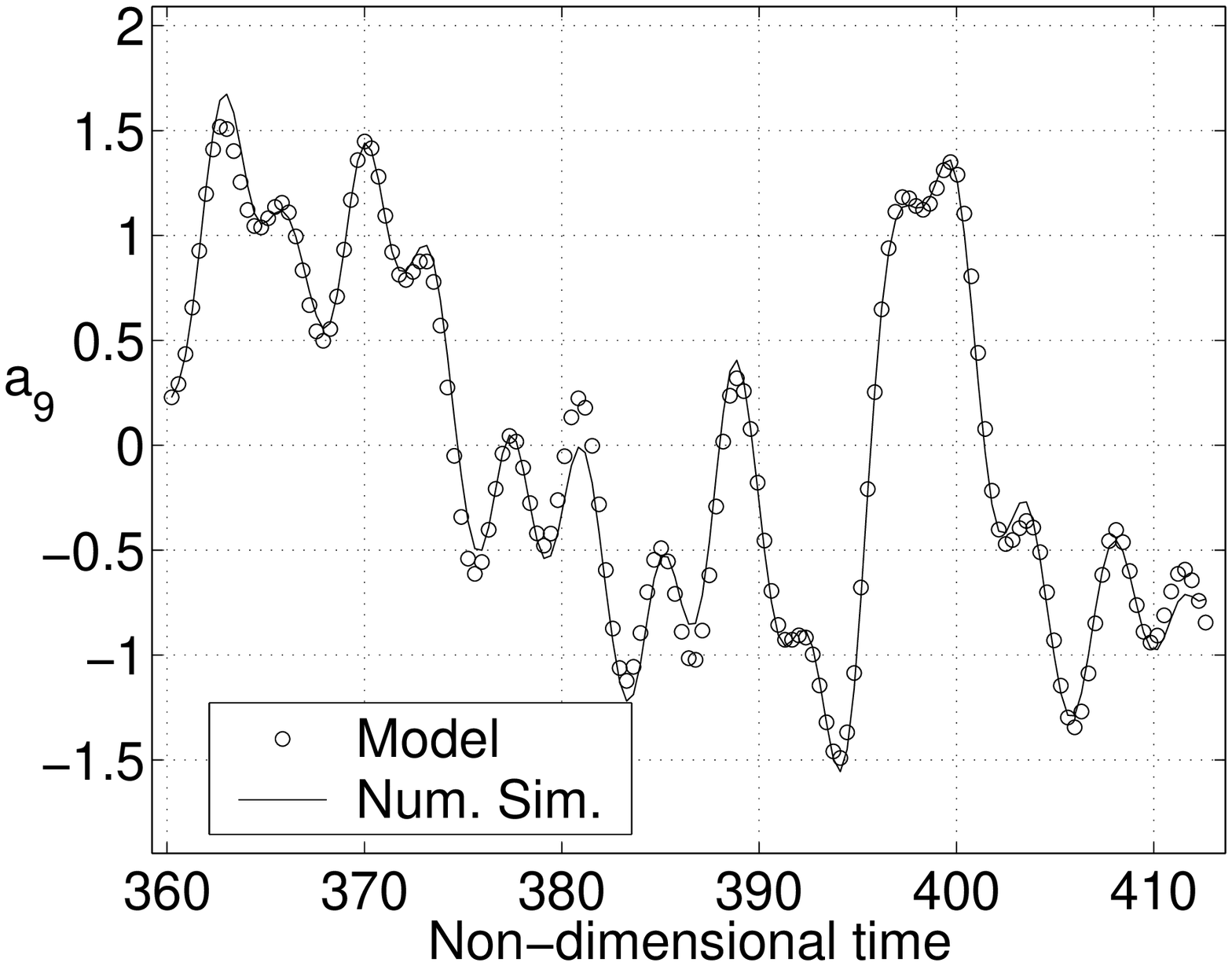}}
  \subfigure{
    \includegraphics[width=3cm,height=3cm]{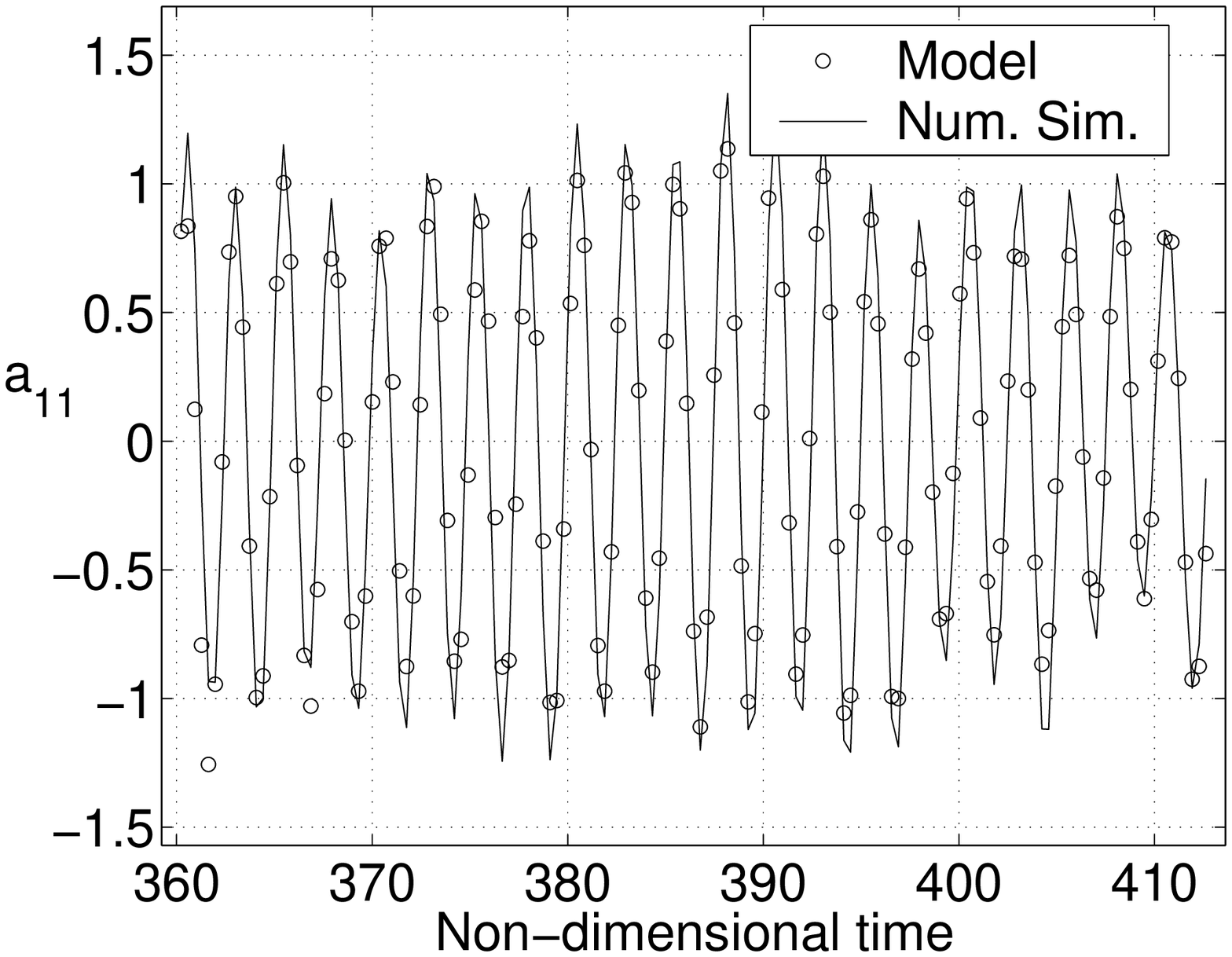}}
  \subfigure{
    \includegraphics[width=3cm,height=3cm]{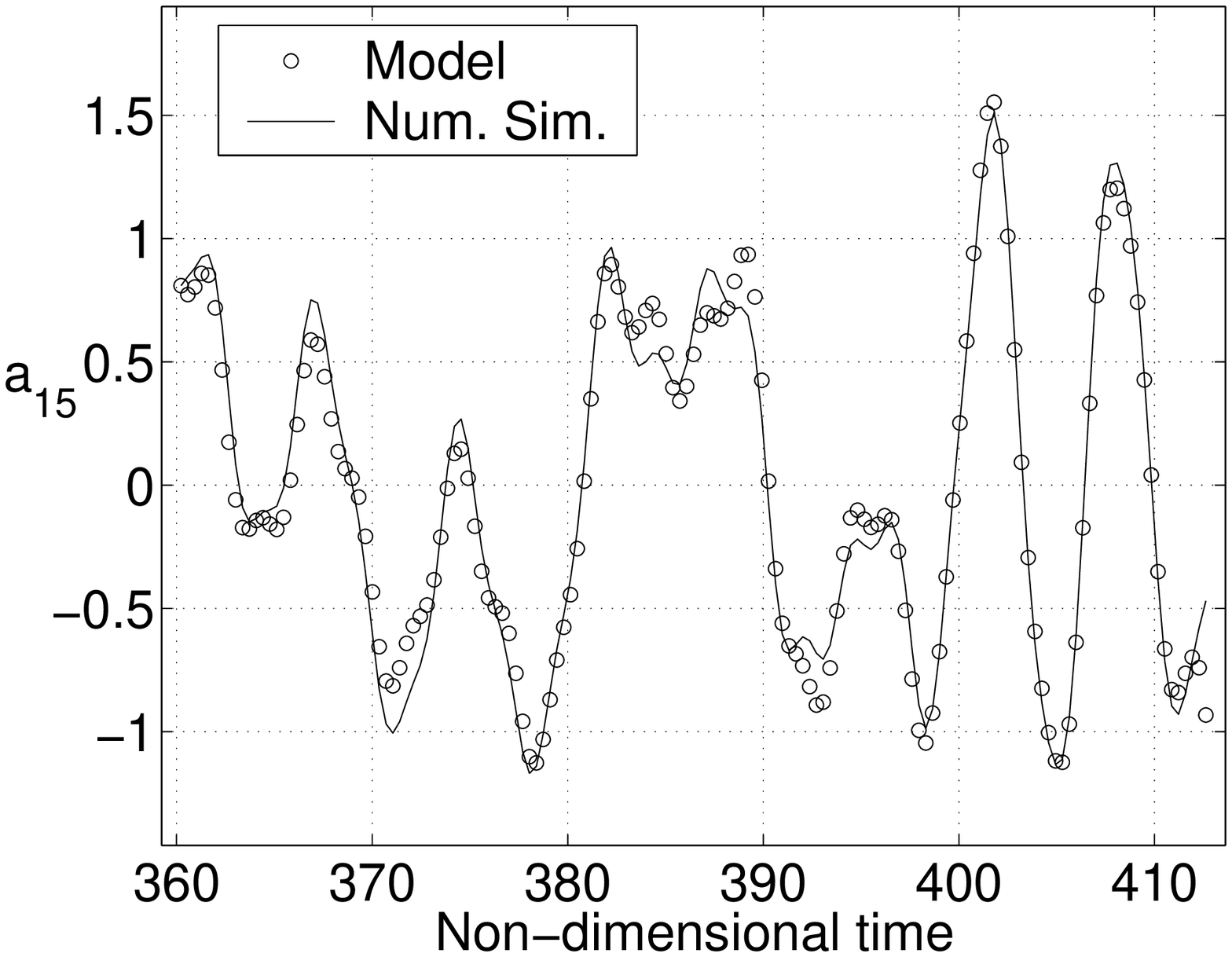}}
  \subfigure{
    \includegraphics[width=3cm,height=3cm]{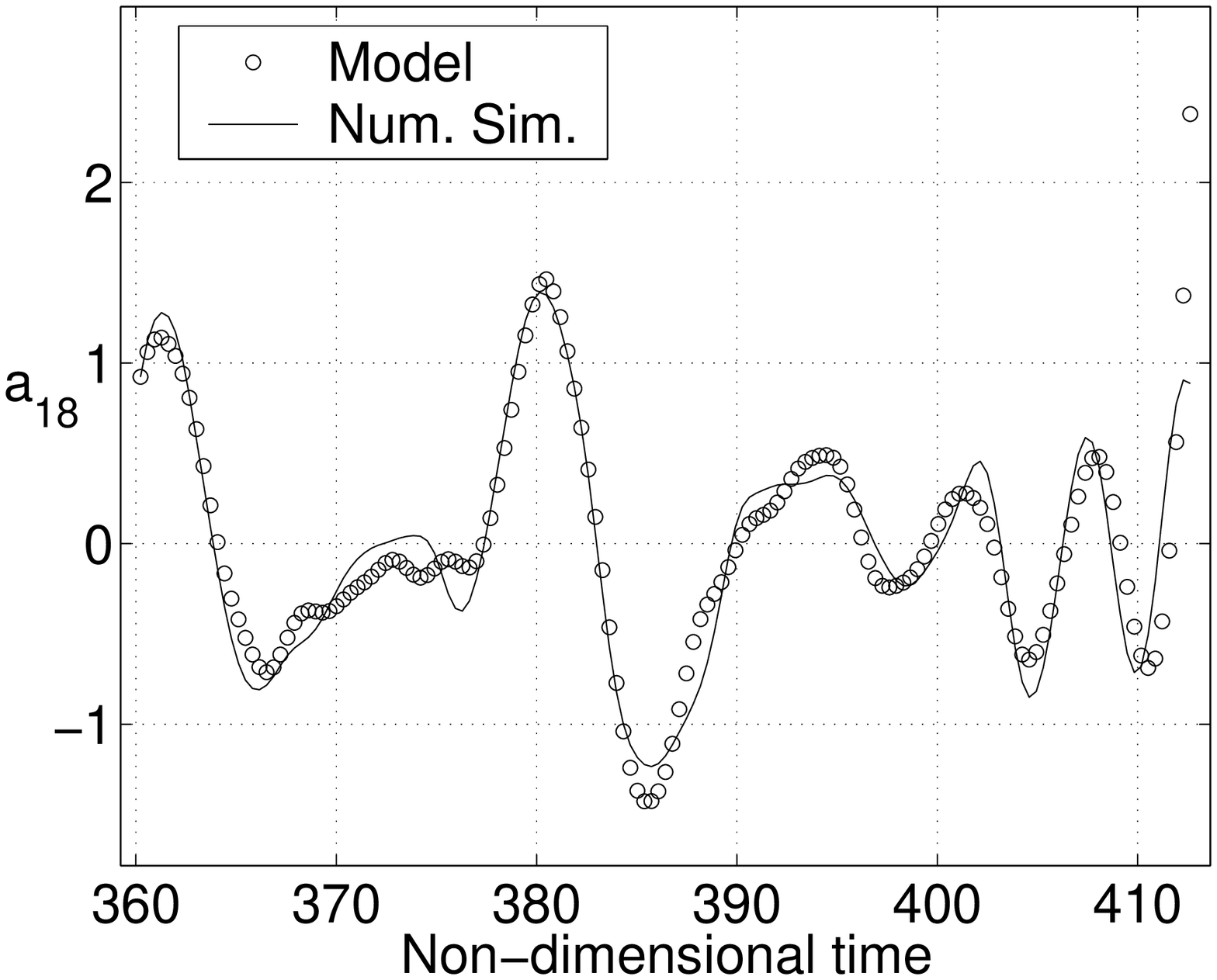}}
  \subfigure{
    \includegraphics[width=3cm,height=3cm]{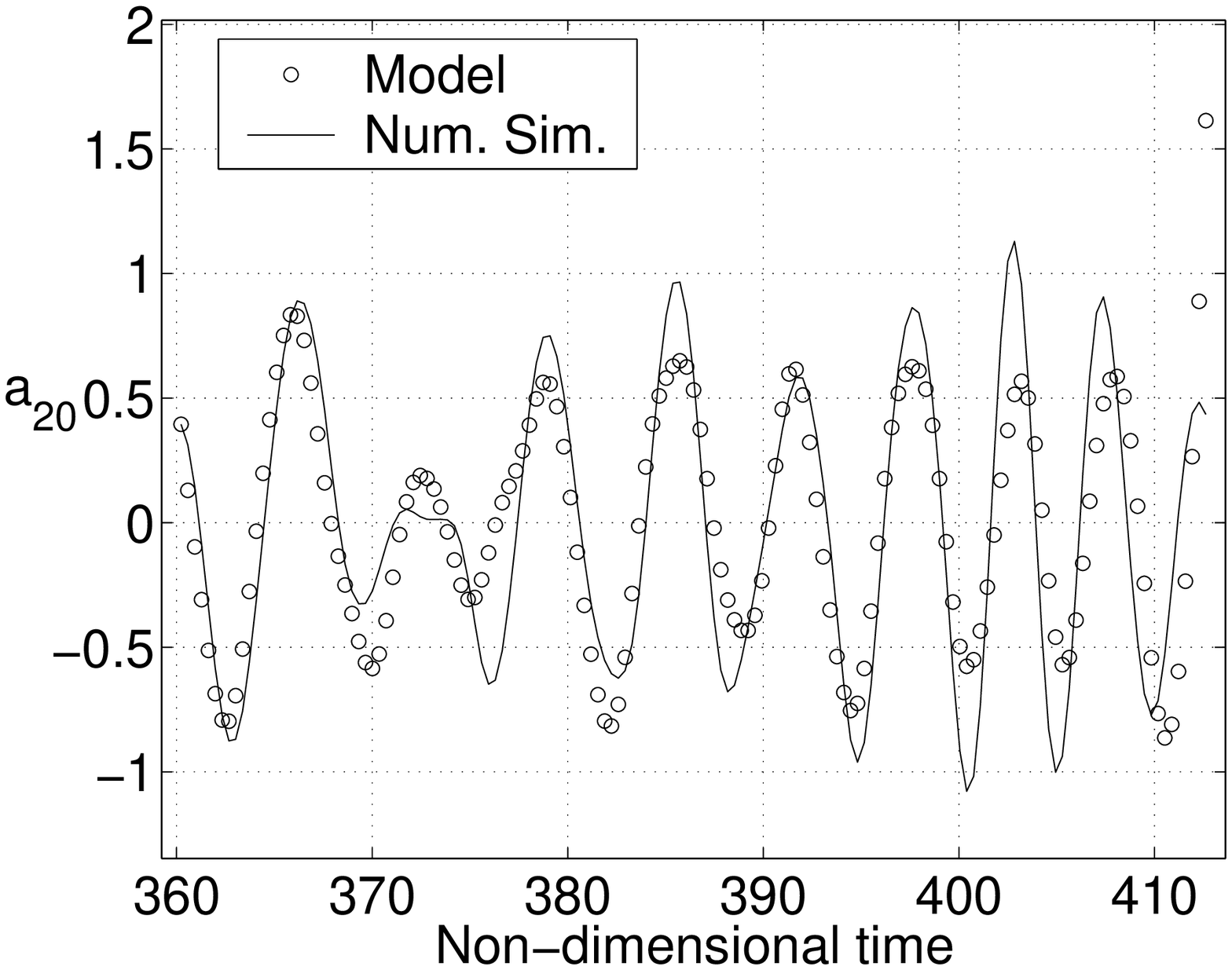}}
\caption{Three-dimensional flow: 
  projection of the fully resolved Navier-Stokes simulations over the POD
  modes (continuous line) vs. the integration of the
  dynamical system inside the calibration interval,
  obtained retaining the first 20 POD modes
  (circles). Only 8
  representative modal coefficients are shown here.} 
\label{fig:calPOD2}
\end{figure}

The resulting model for the two-dimensional flow configuration
considered here is  very accurate in describing the asymptotic
attractor (\cite{Galletti2004} and  \cite{Galletti2006}). For the
three-dimensional case, it was shown in \cite{Buffoni2006}  
that the calibrated model is capable of accurately reproducing the
complicated flow dynamics resulting from the interaction of the
three-dimensional vortex wake with the confining walls inside the
calibration interval. This is shown in figure~\ref{fig:calPOD2}, where
the predictions of some POD modal coefficients given by the dynamic
POD model within the calibration interval are compared to those
obtained from the projection  of the fully resolved Navier-Stokes
simulations. 
However, these results rapidly deteriorate as soon as the model
becomes an actual prediction tool, i.e., outside the time interval in
which the calibration is performed. Nonetheless, the model can be used
in conjunction with experimental measurements in order to reconstruct
the entire flow field, as will be shown in the next section.
  
\section{Non-linear observer}

Our aim is to provide an estimation of the modal coefficients $a_i(t)$
starting from $N_s$ flow measurements $f_k,~k \in
\{1,\dots,N_s\}$. Let $\bar{\alpha}_i(t)$ be the projection of the
velocity field $\vct{u}(t)$ over the $i$-th POD mode and $\alpha_i(t)$
be its estimated value at time $t$. 
 
We assume that each measurement $f_k$ is a scalar quantity which
depends linearly on the instantaneous velocity field $\vct{u}(t)$. For
instance, $f_k$ can be a point-wise measurement of a velocity
component, or of a shear stress, or it can be a spatial average of a
linear combination of velocity components. 

The available spatial information may be exploited by using a LSQ
approach, as done in \cite{Galletti2004}. At any given time $\tau$,
thanks to the linearity of $f_k$ with respect to $\vct{u}$ and to the
modal decomposition of the velocity field (see Eq.~(\ref{galrap})),
$f_k$ can be written in terms of POD modes 
\begin{equation}
f_k\tond{\vct{u}\tond{\tau}} \simeq \sum_{j=1}^{N_r} a_j(\tau)
f_k\tond{\vct{\Phi}^j} 
\end{equation}
where $f_k\tond{\vct{\Phi}^j}$ is obtained from the application of
$f_k$ to the vector field associated to mode $\vct{\Phi}^j$. Then, the
following least-squares problem has to be solved for every $\tau$ 
\begin{equation}
\label{minq}
\min_{\{a_1(\tau),\dots,a_{N_r}(\tau)\}
}
\sum_{k=1}^{N_s} {\left(
f_k\tond{\vct{u}\tond{\tau}} - \sum_{j=1}^{N_r} a_j(\tau)
f_k\tond{\vct{\Phi}^j} 
\right)^2}
\end{equation}
This problem leads to the solution of a $N_r$-dimensional linear
system of equations. Once this problem is solved, the estimated modal
coefficients can be written    
\begin{equation}
\label{eq:lsq}
\alpha_j(\tau)=\sum_{k=1}^{N_s} \Upsilon_{k j} 
f_k\tond{\vct{u}\tond{\tau}}
\end{equation}
where $\Upsilon$ is a known rectangular matrix of size $N_s \times
N_r$. The error minimization (\ref{minq}) leads to a linear
representation of the estimated modes as a function of the
measurements. 

The LSE approach, conversely, exploits temporal rather than spatial
information and is based on the assumption that a linear relation
exists between the modal coefficients and the measurements
\begin{equation}
\label{eq:lse}
\alpha_j(\tau)=\sum_{k=1}^{N_s} \Lambda_{k j} 
f_k\tond{\vct{u}\tond{\tau}}
\end{equation}
where $\Lambda$ is now an unknown rectangular matrix of size $N_s
\times N_r$. This matrix is determined by imposing the condition that
$\forall j \in \{1,\dots,N_r\}$ and $\forall k \in \{1,\dots,N_s\}$ 
\begin{equation}
\label{eq:callse}
\int_0^T {\bar{\alpha}_j(t) \, 
f_k\tond{\vct{u}\tond{t}}
\, dt} = \int_0^T {\sum_{m=1}^{N_s} \Lambda_{m j} 
f_m\tond{\vct{u}\tond{t}}
f_k\tond{\vct{u}\tond{t}}~dt
}
\end{equation}
The time interval $[0,T]$ is the same as that considered for building
the POD modes. Hence, since the left-hand side is known, a set of
linear equations is obtained; these uniquely define the matrix
$\Lambda$.  

The LSQ and and LSE both provide linear estimation of the modal
coefficients. Matrices $\Upsilon$ and $\Lambda$ have the same size,
although the coefficients are different. In the following we overcome
the assumption of a linear relation. 

Let us assume that a certain number of measurements at consecutive
times $\tau_m$, $m\in\{1,\dots,N_m\}$ are available. The main idea of
the dynamic-estimation approach proposed here is to impose that the
coefficients of the modal expansion of the velocity field give the
best approximation to the available measurements, using either LSQ
(\ref{eq:lsq}) or LSE (\ref{eq:lse}), and that  at the same time they
satisfy as closely as possible the non-linear low-order model
(\ref{model}).  

In the LSQ case this is done by  minimizing the
sum of the residuals of (\ref{eq:lsq}) and the
residuals of (\ref{model}) for all times
$\tau_m$.
More precisely, let $\vct{\alpha}(t): \mathbb{R} \rightarrow
\mathbb{R}^{N_r}$ and 
$\vct{\alpha}(t)=\{\alpha_1(t),\dots,\alpha_{N_r}(t)\}$, we have 
\begin{equation}
\label{klsq}
\vct{\alpha}(t)=  \mathop{\mathrm{argmin}}_{\vct{a}(t)} 
\sum_{m=1}^{N_m}{\left(
C_R
\sum_{r=1}^{N_r}R^2_r(\vct{a}(\tau_m))+\sum_{r=1}^{N_r}(
a_r(\tau_m)-\sum_{k=1}^{N_s} \Upsilon_{k r} 
f_k\tond{\vct{u}\tond{\tau_m}} 
)^2 \right)}
\end{equation}
where $\vct{a}(t)=\{a_1(t),\dots,a_{N_r}(t)\}$. 
The parameter $C_R$ 
weights more 
the measurements (LSQ) or  the dynamic model in the definition of the
residual norm. 
It could be systematically tuned, or it could be a matrix.
In the numerical experiments
reported here, this parameter has been
in a heuristic way, leaving a consistent analysis to future investigations.
The minimization of this functional is reduced to a
non-linear algebraic problem. As in \cite{Galletti2006}, a pseudo-spectral
approach is used and  
each $a_r(t)$ is expanded in time using Lagrange polynomials
defined on Chebyshev-Gauss-Lobatto collocation points. 
The necessary conditions for the minimum result in a non-linear set of
algebraic equations for the coefficients of the Lagrange
polynomials. The solution is obtained by a Newton method, which, in
the present applications,  
usually converges in a few (typically 5 to 8)
iterations.
The solution of the problem~(\ref{klsq}) provides
an estimation for the POD modal coefficients for all modes
and for all instants at which measurements are available. 
This allows the reconstruction of the
entire flow field at the same instants
through equation~(\ref{galrap}). The above method, therefore, 
represents a non-linear observer of the flow
state. In the following, it  will be referenced as K-LSQ. 

A similar approach can be obtained 
for the LSE technique, by substituting in Eq.~(\ref{klsq}) the residuals of 
Eq.~(\ref{eq:lse}) instead of those of Eq.~(\ref{eq:lsq}). This approach is
referenced as K-LSE.   

In literature, there exist other flow estimation techniques that are
non-linear in the flow measurements. In the following we will compare
the results of the proposed non-linear dynamic estimation to one of
them, a quadratic extension to LSE (\cite{adr77},
\cite{naguib2001}). This method is based on the assumption that
equation (\ref{eq:lse}) is just the first term of a Taylor expansion
with respect to the sensor measurements, whereas QSE takes into
account the second order term, too. Hence, we have 
\begin{equation}
\label{eq:qse}
\alpha_j(\tau)=\sum_{k=1}^{N_s} \Lambda_{k j} 
f_k\tond{\vct{u}\tond{\tau}}+ \sum_{k=1}^{N_s} \sum_{m=1}^{N_s}
\Omega_{k m j} 
f_k\tond{\vct{u}\tond{\tau}} f_m\tond{\vct{u}\tond{\tau}}
\end{equation}
where the scalar coefficients $\Lambda_{k j}$ and $\Omega_{k m j}$ are
obtained using double, triple and quadruple correlations between
measurements in an equation equivalent to (\ref{eq:callse}). This
approach is referred to as QSE.

Once the matrices appearing in equations (\ref{eq:lsq}) (\ref{eq:lse})
and (\ref{eq:qse}) are computed, the estimation of the modal
coefficient at a certain time is based on the measurements made at the
same time. In contrast, \cite{ewing99}, \cite{tinney} proposed to take
into account integrated temporal data by assuming a linear dependence
between the modal coefficients and the flow measurements in a
non-local way, by working in the frequency domain. Let $\hat{\alpha}$
be the Fourier transform of $\alpha$ and $\hat{f}_j$ that of $f_j$,
then for each frequency we pose 
\begin{equation}
\label{eq:sls}
\hat{\alpha}_j=\sum_{k=1}^{N_s} \hat{\Gamma}_{k j} 
\hat{f}_k
\end{equation}
where $\hat{\Gamma}_{k j}$ is a matrix obtained by appropriate
ensemble averages and depends on the frequency. In the time
domain this amounts to a convolution integral between the measurements
and the time dependent matrix $\Gamma$.  We call this approach SLSE.
As compared to QSE and SLSE, the dynamic estimation procedure that we
propose is non-linear and, at the same time, it takes into
account the evolution of the modal coefficients in time by
constraining such evolution to a model, in the weak sense determined
by \ref{klsq}.

Concerning the applicability of the methods described above, it is
important to recall that the LSE and LSQ approaches are readily
applicable to real-time estimation, as well as the QSE, although the
cost of this last approach scales as $N_s^2$ instead of linearly as in
the previous two cases.
Conversely, the SLSE approach is more difficult to be used for
real-time estimations, since it uses the whole temporal 
history of the measurements,
collected in a time interval,
coupled together (linearly) 
via the Discrete Fourier Transform (DFT). 
This implies that 
the estimation problem must be tackled after having collected enough
temporal information
and it consists of a number of LSE problems equal to the
number of retained frequencies, plus additional DFT's of the
measurements and of the estimated POD coefficients. 

Similarly to what done in the SLSE approach, in the present 
dynamic estimations
the temporal histories of the measurements are
coupled together (non-linearly) by the dynamic POD model. 
This aspect poses difficulties in a real-time application. Indeed,
as pointed out in the introduction, K-LSQ or K-LSE
are thought to be applied a-posteriori,
because their computational cost, 
although unimportant in a post-processing phase, is 
large for a real-time analysis.
Nevertheless, although actual real-time applications are premature, 
a proposal for their prospective implementation for 
real-time estimation is the following. The flow state at a given
time $t^*$ could be estimated by considering 
the measurements taken at that time and at the previous $N_m-1$ ones. 
At the successive sampling time, the corresponding
new measurements are added and the oldest
ones are dropped, keeping the number of measurements considered constant.
In other words, reconstruction is carried out using a fixed
number of measurements distributed in a time interval which is 
located before $t^*$, and which translates as time increases. 
The sampling rate (i.e. $\tau_m-\tau_{m-1}$)
and $N_m$ can be tuned in order to decrease the
computational costs while granting the level of accuracy required by
the particular application. Moreover, when a new set of measurements
is added, the Newton method for solving the non-linear system would
be restarted from the previous solution, which is already close to the
final solution, thus definitely reducing the number of iterations for
convergence. 

In contrast with the other methods, the proposed approaches
need a working Galerkin model as a fundamental ingredient. 
The construction of such a model can be carried out from
the information needed to build the POD database, a necessary step
for all the methods considered here. Therefore no additional
information is needed as compared to other approaches. 

\section{Results and discussion}

The K-LSQ and K-LSE are used to
reconstruct the flow in the configuration described in
Sec.~\ref{sec:flow}, both in the two- ($\Rey=150$) and
three-dimensional ($\Rey=300$) cases. Results are
compared to those obtained by  
the most common techniques available in the literature  
and cited in the Introduction.

Accuracy in the prediction of the single modal coefficients
and in the reconstruction of the velocity fields were appraised.
In both cases, differences with respect to the reference case (DNS)
were quantified in
terms of relative error in the $L^2$ norm, i.e., the
$L^2$ norm of the difference between the estimated and the reference
quantity divided by the norm of the reference quantity.

Several parameters are involved in  the set-up of the K-LSQ and K-LSE
models. They are related to (i) the dynamic POD model: number of
retained modes, 
calibration interval, number and temporal distribution of available
snapshots; (ii)  the selected flow measurements: number, type and
collocation.

The parameters involved in the derivation of the POD dynamical system
used herein are the same that were selected in previous studies 
(\cite{Galletti2004}, \cite{Buffoni2006}). 
As for the flow measurements, both velocity and shear-stress sensors
were used. 
While velocity measurements are often considered in the literature,
due to their widespread use in practice, 
shear-stress sensors are less
common. Nevertheless, they were used here mainly because they are 
challenging from a numerical point of view, as
they involve spatial derivatives of the POD modes. Also, they can be
implemented in practice although limitations of accuracy and time
resolution may exist (see, for instance, \cite{spazza}).
Different sensor locations were tested, 
to account for the sensitivity of the
proposed approaches to sensor placement.   
Since 
the performance of the standard techniques such as LSE or LSQ is
influenced by sensor placement, 
some sensor configurations were selected 
following the suggestions given for LSE in \cite{Cohen2004}. 
On the other hand, none of the considered sensor
configurations is optimized for K-LSE or K-LSQ, in order to verify the 
sensitivity of such methods  with respect to sensor placement. In fact,
optimal sensor placement may turn out to be a time-consuming
operation for complex three-dimensional flows.

\subsection{Two-dimensional case: $Re=150$}

For this rather simple flow, we consider the situation in which a
limited number of measurements are available,
i.e. only 2 sensors. 
Three different configurations were analyzed, two involving
streamwise velocity sensors  and one involving shear-stress sensors.

The velocity sensors were placed 
in relation to the spatial structure of the streamwise component of the
first two POD modes. 
In particular, in the first configuration one streamwise velocity
sensor is placed on the 
maximum of the first POD mode which is closest to the cylinder 
($P1\simeq (2.39,0.52)$) and 
one in the middle between P1 and the minimum of the second POD mode
closest to the cylinder ($\simeq (1.96,0.50)$). 
The second configuration has the first streamwise velocity sensor in
P1 and the second one 
in the point $(1.98,-0.76)$.
A third configuration was considered with two
shear-stress sensors located 
on the confining walls ($y = \pm 4.0$) at $x = 4$, in a region which
satisfies the following  criteria on a shedding cycle:   
the rms value of
the shear-stresses is maximum and  the
reconstruction error of the shear-stresses is minimum for a given
number of POD modes. 

The POD low-order model of the two-dimensional flow is obtained
by using 95 snapshots, uniformly distributed throughout
two vortex 
shedding cycles ($T \simeq 13$ is the non-dimensional duration of the
time interval), and by retaining $N_r=6$ modes. The calibration of the
model is performed
in the same interval using 81 collocation points.   
As shown in \cite{Galletti2004}, the calibrated model 
accurately reproduces the flow inside and outside the
calibration interval.  

\begin{figure}[h!]
  \centering
    \subfigure[]{\includegraphics[width=6.5cm]{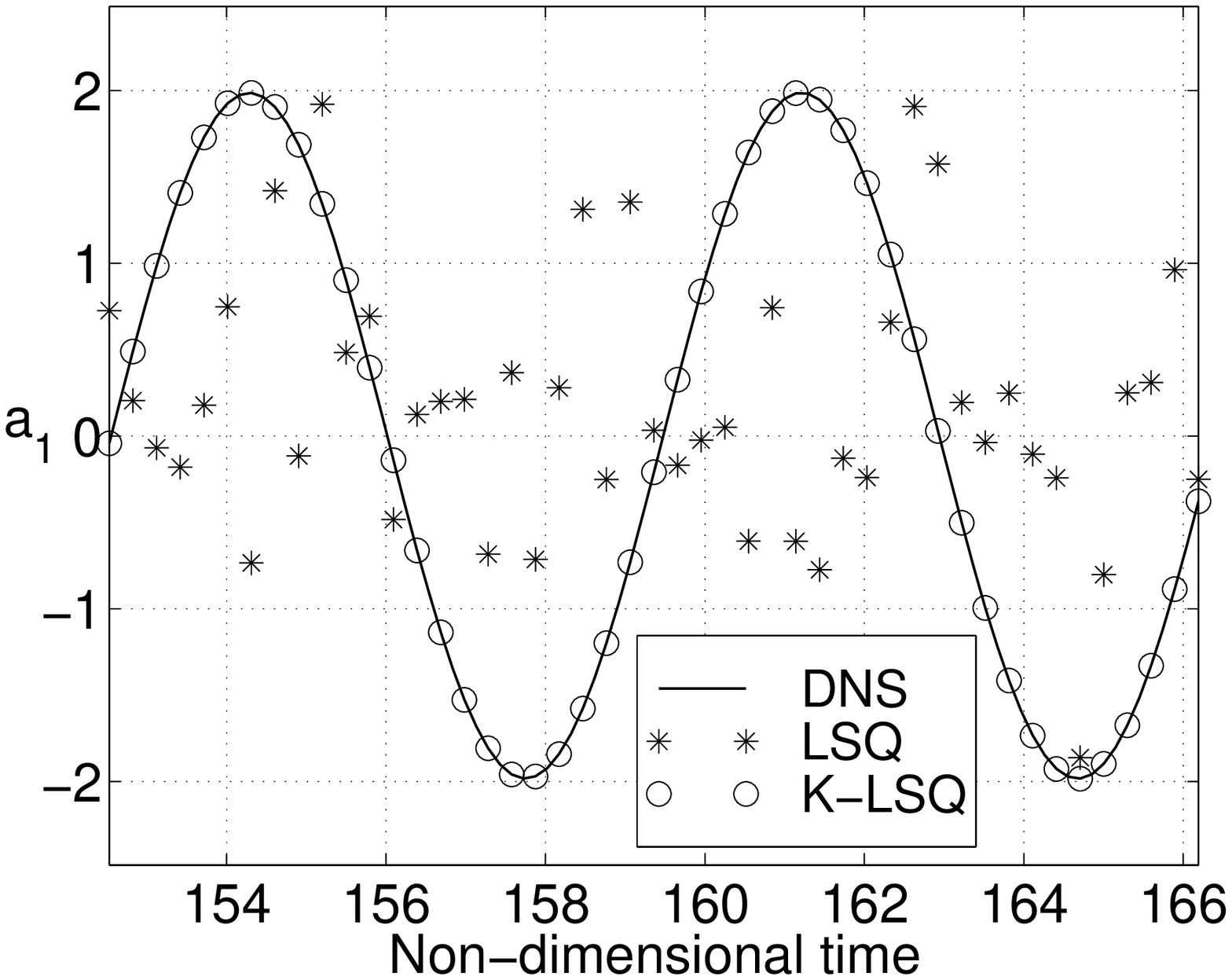}}%
    \hspace{4.mm}
    \subfigure[]{\includegraphics[width=6.5cm]{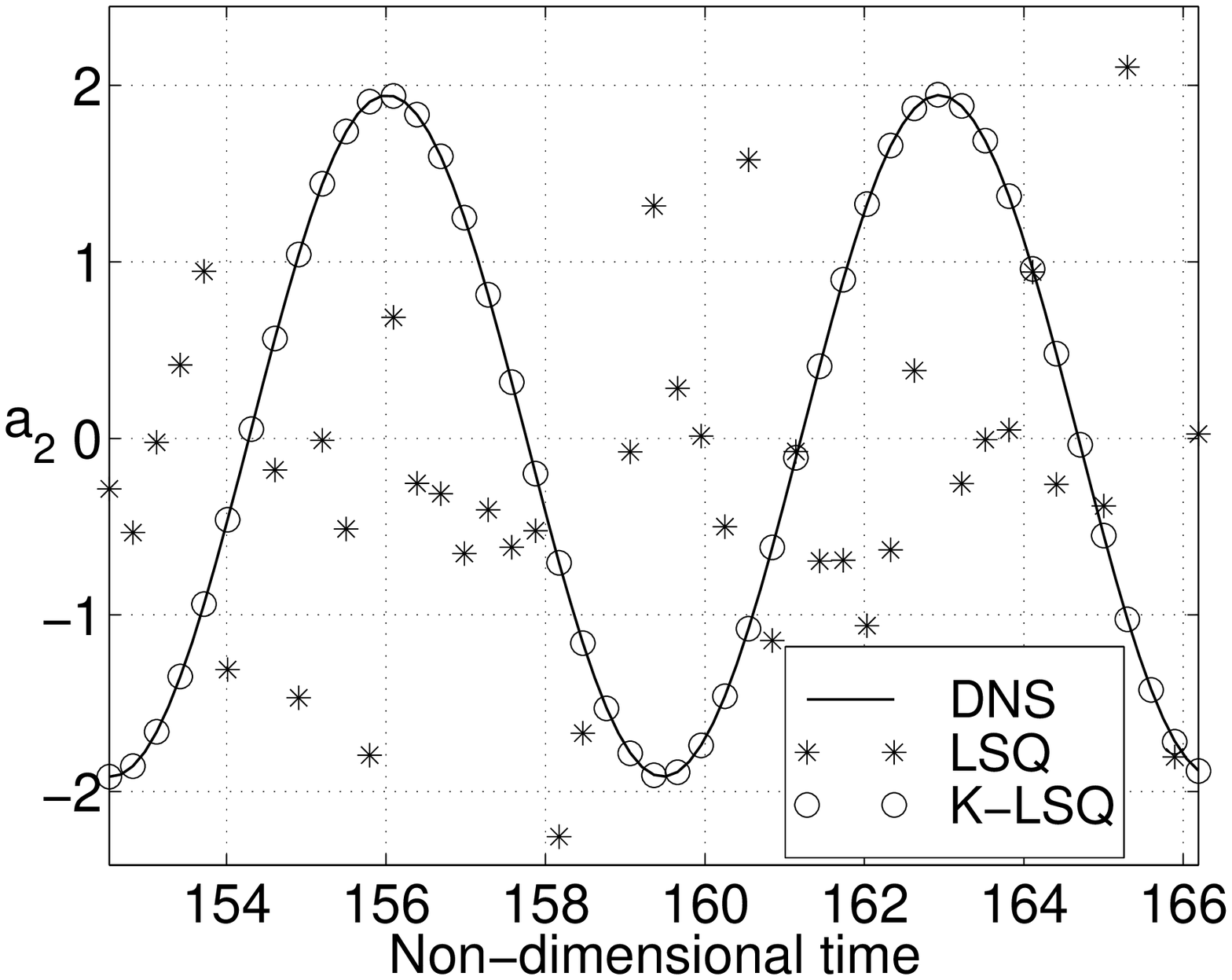}}\\[2.mm]%
    \subfigure[]{\includegraphics[width=6.5cm]{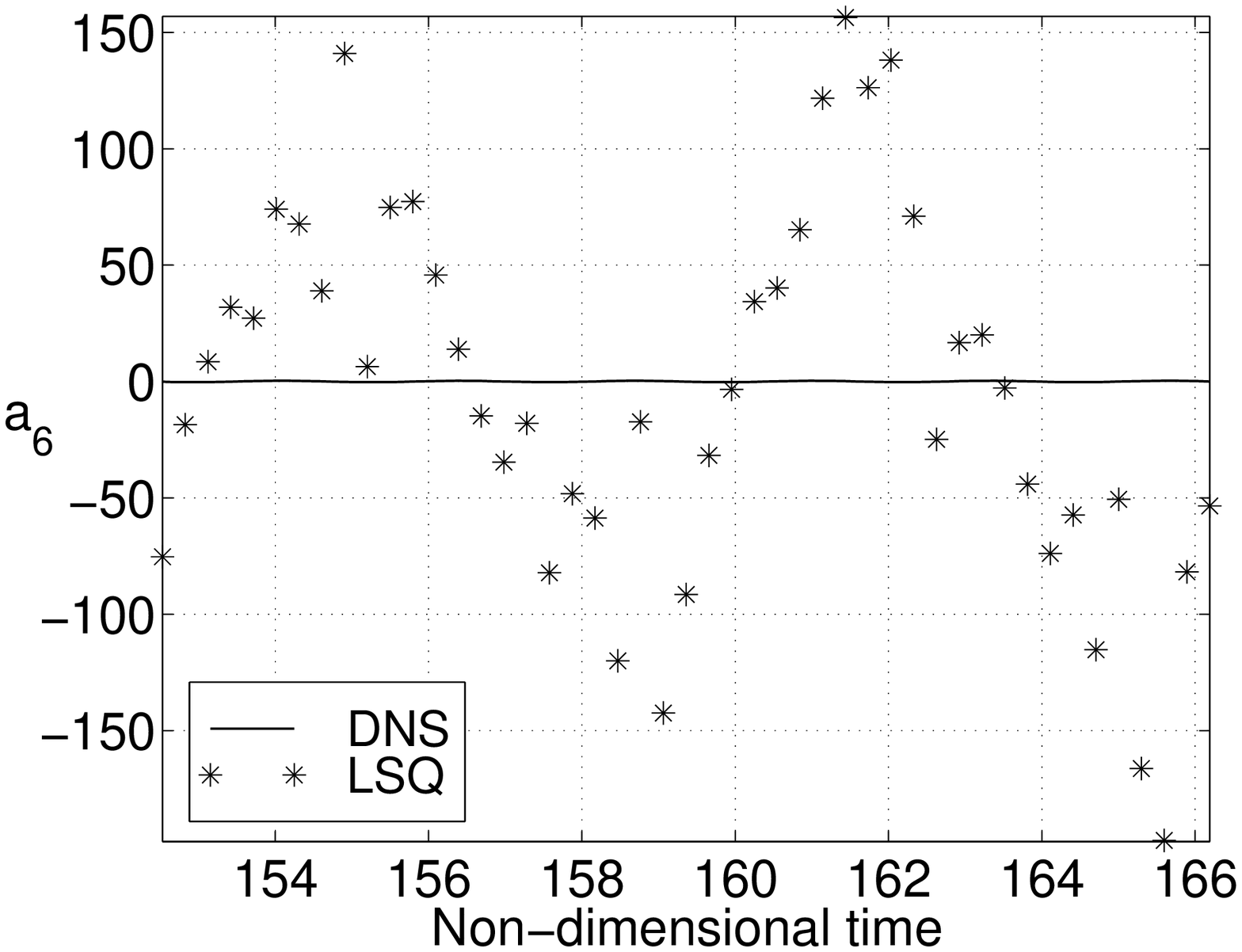}}%
    \hspace{4.mm}
    \subfigure[]{\includegraphics[width=6.5cm]{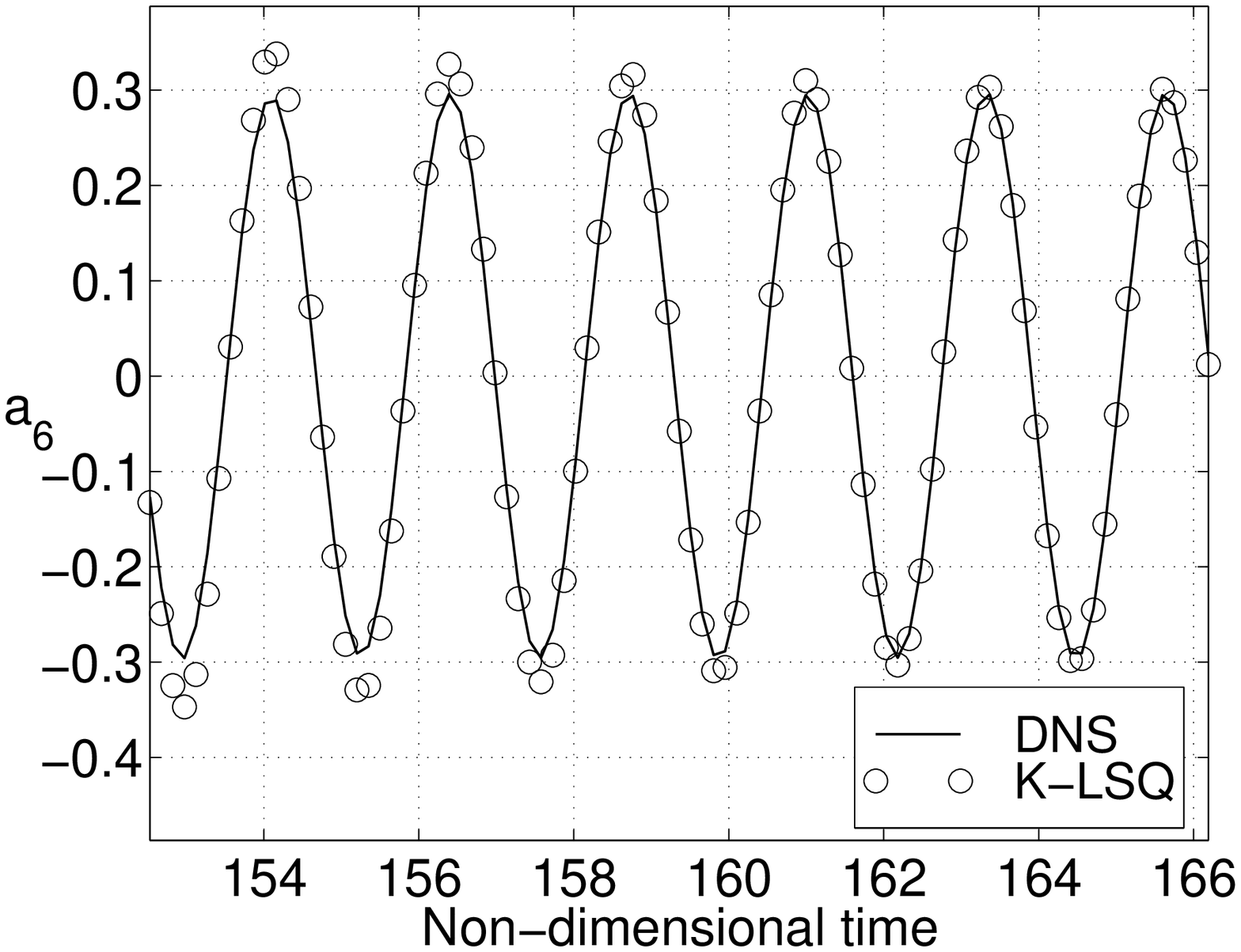}}%
  \caption{\em (a) and (b): POD modal coefficients $a_1$ and $a_2$ estimated
    with the K-LSQ and the LSQ approaches, compared to their exact
    values computed by projection of the DNS velocity fields. 
    (c) coefficient $a_6$, DNS and LSQ; (d) coefficient $a_6$, DNS and
    K-LSQ. Note that in (c) and (d) different axis scales are used.}
  \label{2d}
\end{figure}

Figure~\ref{2d} shows the POD modal coefficients predicted by K-LSQ
and by LSQ considering the third (shear-stress sensors) configuration,
together with their reference values found from DNS. In addition,
errors in the prediction of the modal coefficients given by 
LSQ, LSE, QSE, K-LSQ 
and K-LSE in the first (velocity sensors) and third 
(shear-stress sensors) configurations are reported in
Table~\ref{res2dcoeff}(a) and (b), respectively. 
The values obtained for the second considered sensor
configuration are not shown since they are very similar to those of
the first one.

The time interval over which
reconstruction is performed is approximately 13 time-units long
(non-dimensional time); it contains two shedding 
cycles, and it starts 
just after the end of the time interval on which the POD model was
built and calibrated. 

\begin{table}
\centering
\begin{tabular}{|c|c|c|c|c|c|c|}
\hline
 Tab. \ref{res2dcoeff}(a) & $e(a_1)\%$ & $e(a_2)\%$ & $e(a_3)\%$ &
   $e(a_4)\%$  & $e(a_5)\%$ & $e(a_6)\%$ \\%
\hline
LSQ & 108.41 & 340.48 & 538.59 & 2210 & 10900 & 8340 \\%
LSE & 71.06 & 27.12 & 99.71  & 97.97 & 99.91 & 99.93 \\%
K-LSQ & 0.47 & 0.55 & 2.58 & 2.66 & 4.65 & 4.67 \\%
K-LSE & 0.82 & 0.76 & 9.82 & 9.82 & 14.98 &  15.59 \\%
\hline
\end{tabular}

\begin{tabular}{|c|c|c|c|c|c|c|}
\hline
 Tab. \ref{res2dcoeff}(b) & $e(a_1)\%$ & $e(a_2)\%$ & $e(a_3)\%$ &
 $e(a_4)\%$  & $e(a_5)\%$ & $e(a_6)\%$ \\%
\hline
LSQ & 100.33 & 140.60 & 752.10 & 1040 & 5490 & 37800 \\%
LSE & 46.87 & 91.97 & 102.83 & 100.99 & 101.40 & 100.21 \\%
K-LSQ & 0.06 & 0.09 & 6.06 & 6.09 & 9.72 &  9.56 \\%
K-LSE & 2.99 & 3.08 & 7.07 &  8.31 & 17.99 & 18.47 \\%
\hline
\end{tabular}
\caption{Relative percentage errors (in $L^2$ norm) 
  on the estimation of the POD modal coefficients ($e(a_i)$) in the
  first (a) and third (b) sensor configuration. In this case
  time-averaging is carried out over the estimation time period.} 
\label{res2dcoeff}
\end{table}

Tables~\ref{res2dcoeff}(a) and (b) show the relative
reconstruction errors on the velocity components and on 
their fluctuating part.
It appears that two (velocity or shear-stress) sensors 
are not sufficient to
obtain reliable predictions of the modal coefficients by LSQ, LSE,
even if LSE leads to a better estimation than LSQ. Accuracy problems
persist also with the QSE approach, even if in this case
the predictions 
are  more accurate than those obtained with LSE. 
The results obtained by LSE and QSE
are more accurate 
in the configuration with velocity sensors because their
positions were selected to be 
effective for LSE, as already pointed out. 
Nevertheless,
the errors on the estimation of the first two modal
coefficients are still large. 
This leads to severe errors in
the estimation of the fluctuating part of the 
velocity field since the first two
POD modes represent about $94.8\%$ of the fluctuating energy. 
Even if the mean flow energy is important with respect to the fluctuating energy,
errors in the modal coefficients lead to detectable errors in
the reconstruction of the velocity components. 
Note that the 
reconstruction errors
on the vertical component are larger than those on the streamwise one
simply because the contribution
of the mean flow on that component is much lower than in the
streamwise component.
Table~\ref{res2dcoeff} shows that both 
K-LSQ and K-LSE  give an accurate estimation 
not only of the first two modal coefficients, but of all the retained
modes. This leads to a precise estimation of the velocity field as
well as of its fluctuating part. 
Moreover, the accuracy of the results is very similar with
shear-stress and velocity sensors, indicating a weak sensitivity of the
approach with respect to the type and location of the sensors. 
This is not the case for the LSQ, LSE and QSE methods, which show a higher
sensitivity to this aspect, confirming what has already been reported in
the literature.  

\begin{table}
\centering
\begin{tabular}{|c|c|c|c|c|}
\hline
Tab. \ref{res2dener}(a) & $\overline{e(U)}\%$ & $\overline{e(V)}\%$ &
     $\overline{e({U'})}\%$ & $\overline{e({V'})}\%$ \\%
\hline
LSQ  & 106.79 & 937.42 & 906.96 & 1370 \\%
LSE  & 6.32 & 37.14 & 53.96 &  54.26 \\%
K-LSQ & 0.63 & 3.97 & 5.39 & 5.80 \\%
K-LSE & 0.69 & 4.42 & 5.93 & 6.46 \\%
\hline
\end{tabular}

\begin{tabular}{|c|c|c|c|c|}
\hline
Tab. \ref{res2dener}(b) & $\overline{e(U)}\%$ & $\overline{e(V)}\%$ &
     $\overline{e({U'})}\%$ & $\overline{e({V'})}\%$ \\%
\hline
LSQ & 273.22 & 2490 & 2350 & 3640 \\%
LSE & 8.05 & 46.14 & 68.32 & 67.44 \\%
K-LSQ & 0.65 & 4.10 & 5.54 & 6.00 \\%
K-LSE & 0.77 & 4.91 & 6.54 & 7.17 \\%
\hline
\end{tabular}
\caption{Relative percentage errors (in $L^2$ norm) on the estimation
  of the velocity components ($\overline{e(U)}$,$\overline{e(V)}$) and
  of their fluctuating part ($\overline{e(U')}$,$\overline{e(V')}$),
  in the first (a) and third (b) sensor configuration. In this case
  time-averaging is carried out over the estimation time period.}
\label{res2dener}
\end{table}

We compare our results to those of
\cite{Cohen2004},  Tab.~2(a) 13-th case. 
With LSE and 2 sensors they found $e(a_1)\simeq 76.6\%$ and $e(a_2)\simeq
15.1\%$, errors that are similar to those reported in
Table~\ref{res2dcoeff}(a) for LSE. Using the dynamic estimation, the
errors on the same coefficients are two orders of magnitude lower.
Furthermore, using K-LSQ method and two shear-stress sensors
(Table~\ref{res2dcoeff}(b)) the first two modal
coefficients are estimated with an error lower than $0.1\%$, i.e.,
three orders of magnitude lower than LSE. 

\subsection{Three-dimensional case: $Re=300$}

The flow patterns in this case are definitely more complex
than those of the previous one. 
For this reason, 24 flow measurements were used for the
reconstruction procedure, organized in five different
configurations, two involving only velocity measurements and three
involving both velocity and shear-stress measurements. 
In the last three configurations, 
the shear-stress sensors were selected following the
same criterion adopted in the 2D case, they are 14 in number and in
all considered cases   
are symmetrically placed on both the confining walls ($y=\pm 4$) at $x=4$
and $z=\{1.2,1.5,2.7,3,3.3,4.5,4.8\}$.
The placement of the velocity sensors has again been chosen on the
basis of
the spatial structure of the streamwise velocity of the
first 12 POD modes.
The different configurations are listed below,
together with a brief description of the rationale
for the placement of the velocity sensors:
\begin{enumerate}
\item 24 velocity sensors distributed on 6 equispaced slices in the
  axial ($z$) direction; on each slice, the sensors are on the lines
  connecting the maximum and minimum closest to the cylinder 
  of the first two POD modes. On each segment, the
  sensors are approximately in the middle, but
  slightly closer to the extrema 
  of the first POD mode.
\item 24 velocity sensors distributed on 4 equispaced slices in the
  axial ($z$) direction; on each slice, 3 points are selected in the
  region of overlapping between the maxima and minima of the
  low-frequency POD modes (modes 3, 4, 7, 8, 9 and 10) and 3 on the
  overlapping region of the extrema 
  of the vortex shedding modes (modes 1, 2, 5, 6, 7,
  11 and 12) (see \cite{Buffoni2006} 
  for details on the separation between low-frequency and
  vortex-shedding POD modes).
\item 14 shear-stress sensors and 10 velocity sensors distributed on 5
  equispaced slices in the axial ($z$) direction; on each slice, the
  velocity sensors are placed on the maximum and minimum closest to the
  cylinder of the first POD mode.
\item 14 shear-stress sensors and 10 velocity sensors. 6 equispaced
  slices in the axial ($z$) direction are considered. On 4 slices, 2
  velocity sensors are placed as in the previous case. Two sensors are
  placed on the remaining slices, corresponding respectively 
  to the maximum and minimum of the third POD
  mode (low frequency mode). 
\item 14 shear-stress sensors and 10 velocity sensors located in the
  wake, at the points reported in table~\ref{tab:sensdet}.
\end{enumerate}
\begin{table}
\begin{center}
\begin{tabular}{cccc|cccc}
Velocity sensor & x & y & z & Velocity sensor & x & y & z\\ 
1 & 5.01 & 1.03 & 2.00 & 
6 & 5.01 & -1.03 & 4.00 \\ 
2 & 6.99 & 1.04 & 2.00 &
7 & 6.99 & -1.04 & 4.00 \\
3  & 6.03 & -1.04 & 2.00 &
8  & 6.03 & 1.04 & 4.00 \\ 
4  & 6.01 & -0.99 & 2.00 &
9  & 6.01 & 0.99 & 4.00 \\
5  & 5.96 & -0.97 & 2.00 &
10 & 5.96 &  0.97 & 3.99 
\end{tabular}
\end{center}
\caption{Positions of the velocity sensors in the three-dimensional case} 
\label{tab:sensdet}
\end{table}
The low-order model of the developed three-dimensional flow was
derived retaining the first 20 POD modes obtained from a database
of 151 snapshots, uniformly distributed over eight vortex shedding
cycles ($\simeq 52$ non-dimensional time units). Calibration was
carried out in the same time interval,
and the 
results obtained integrating the dynamic model 
within the calibration interval are reported
in figure~\ref{fig:calPOD2}. 
More details of this model can be found in \cite{Buffoni2006} 
(the model is denoted there as POD2).
In contrast to the two-dimensional case,  
the POD model is inaccurate outside the
calibration interval, so it
cannot be used as a predictive tool by itself. 
The estimation of the flow is carried out 
here in two different time intervals, 
one starting just after the end of the time interval 
in which the POD model was calibrated
and the other one being about 52 time units distant from the
calibration interval.
Both the intervals are approximately 30 time units long, including
approximately 4 shedding cycles.
\begin{figure}
  \centering
    \includegraphics[width=10cm]{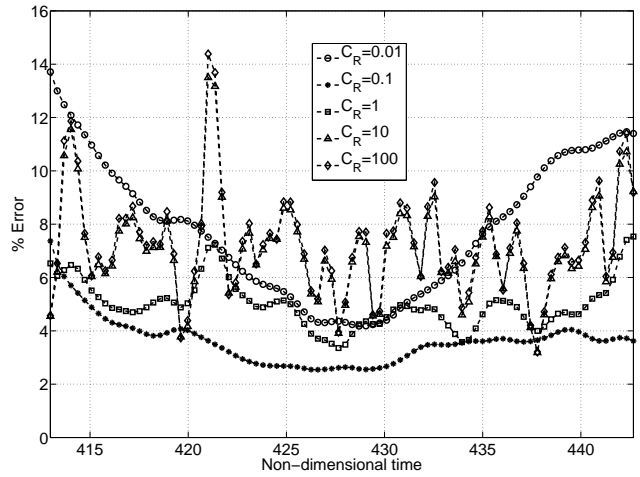}%
  \caption{Relative percent error in the
     reconstruction of the U component
      projected on the retained POD modes as a function of time, when
      varying $C_R$.}
  \label{fig:CR}
\end{figure}
\begin{figure}
  \centering
  \includegraphics[width=5.cm,height=4.25cm]{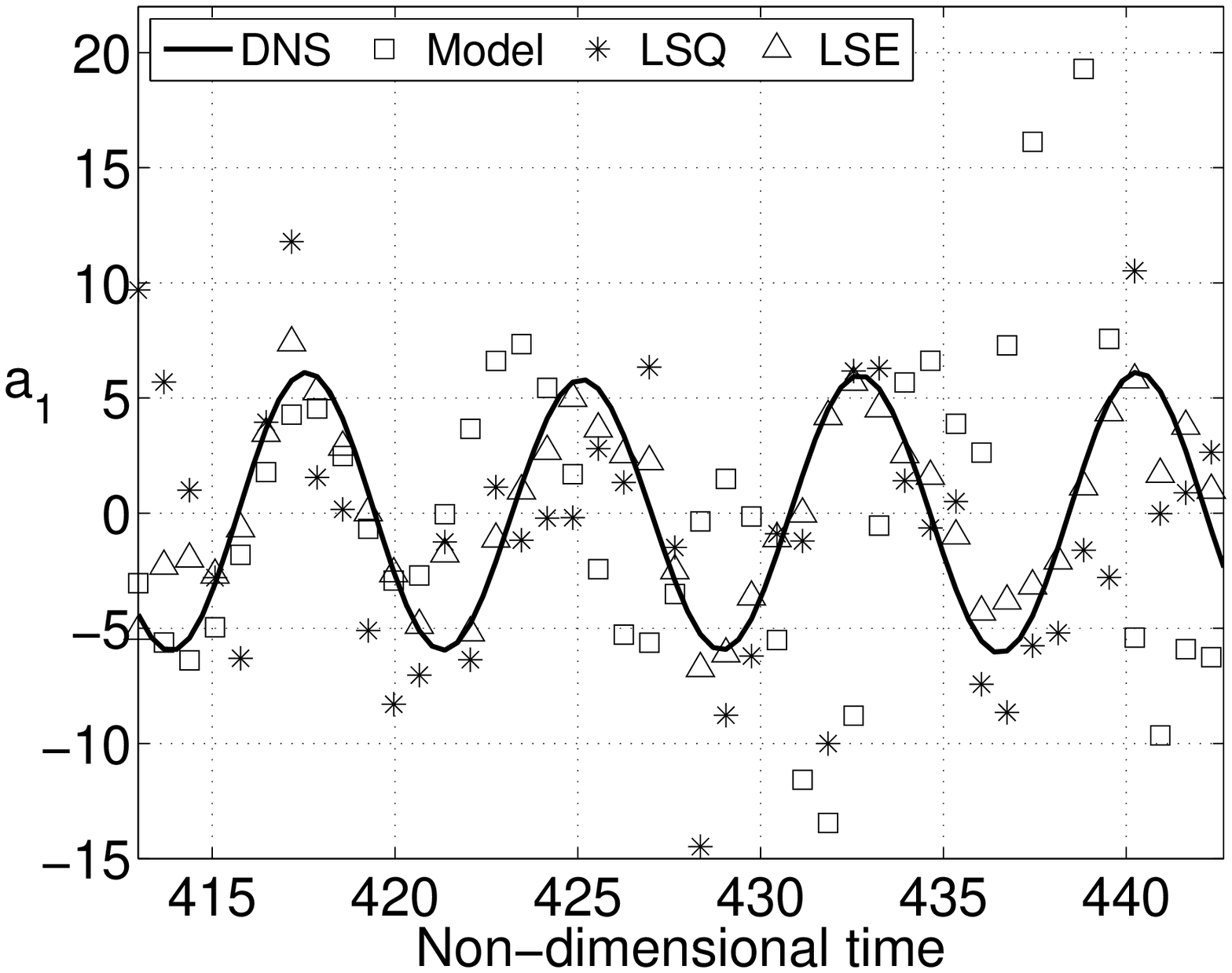}%
  \hspace{1.5cm}
  \includegraphics[width=5.cm,height=4.25cm]{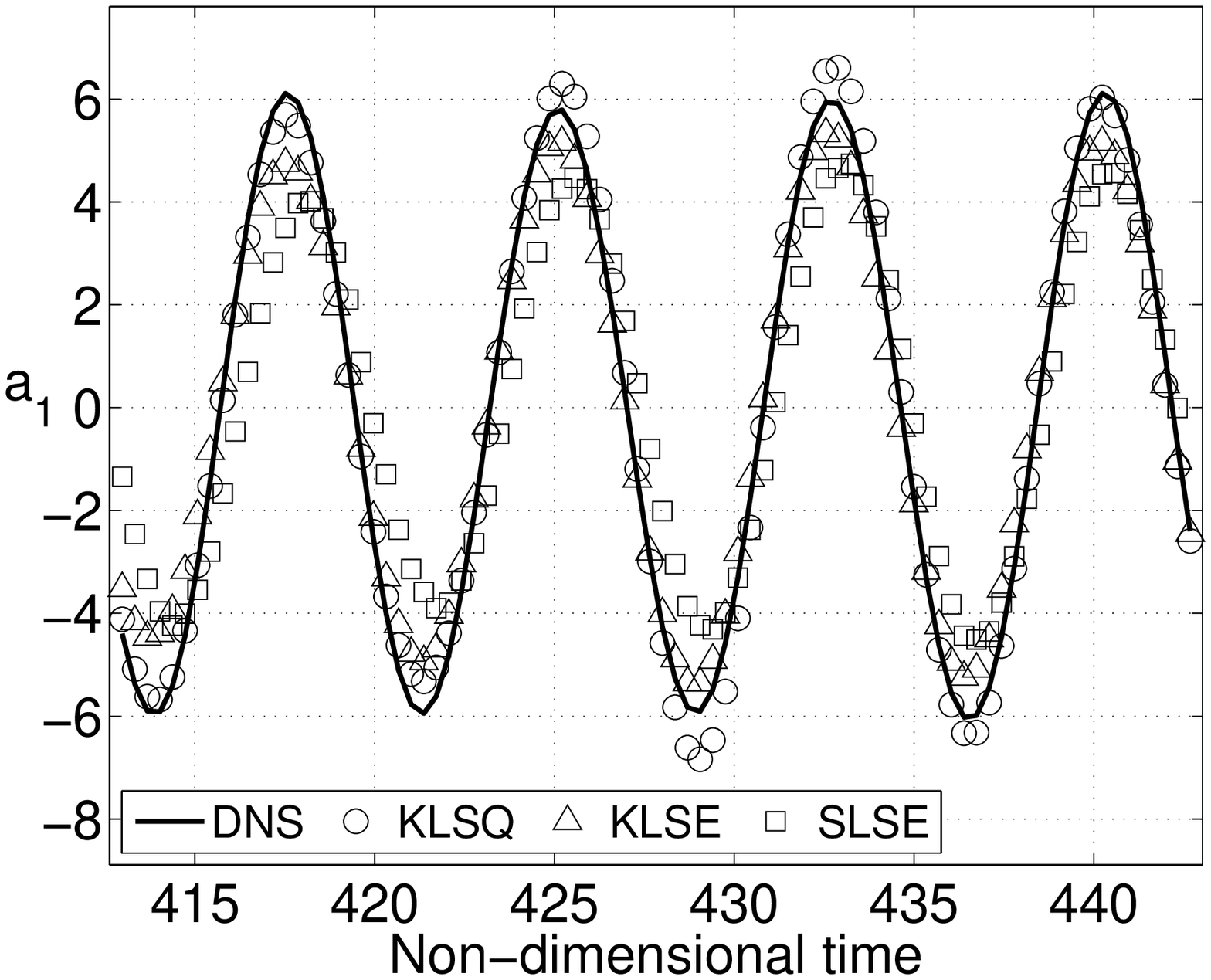}\\%
  \includegraphics[width=5.cm,height=4.25cm]{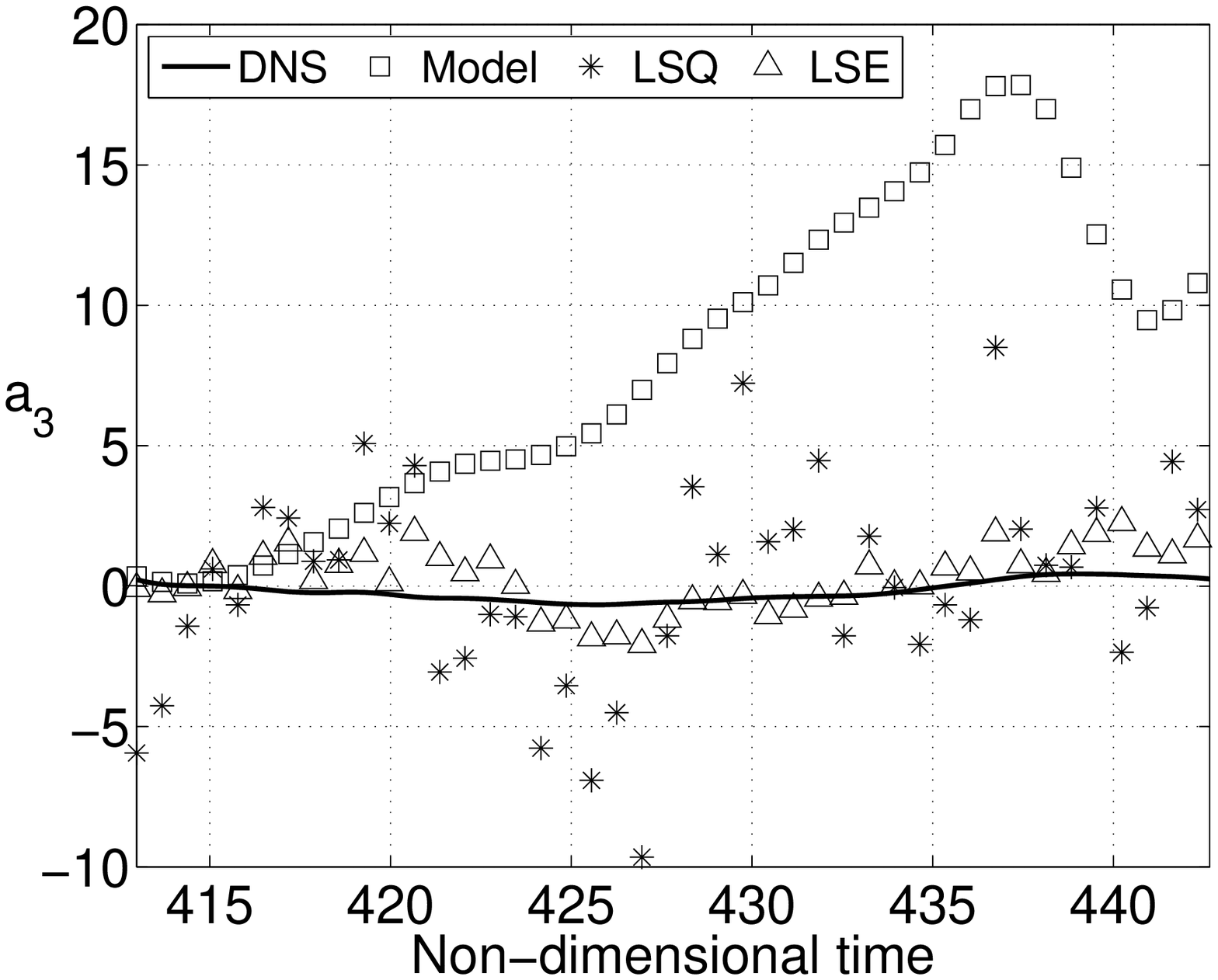}%
  \hspace{1.5cm}
  \includegraphics[width=5.cm,height=4.25cm]{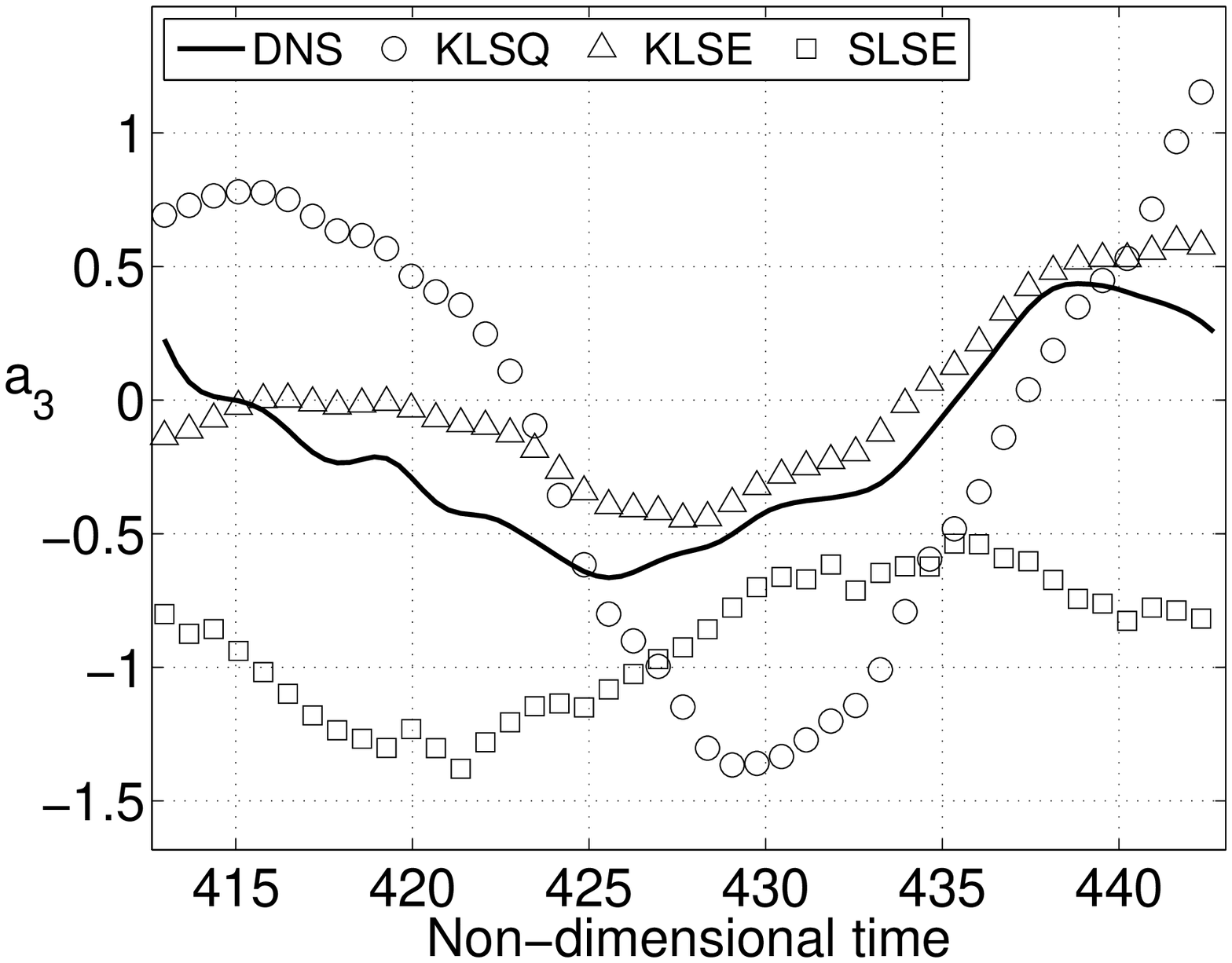}\\%
  \includegraphics[width=5.cm,height=4.25cm]{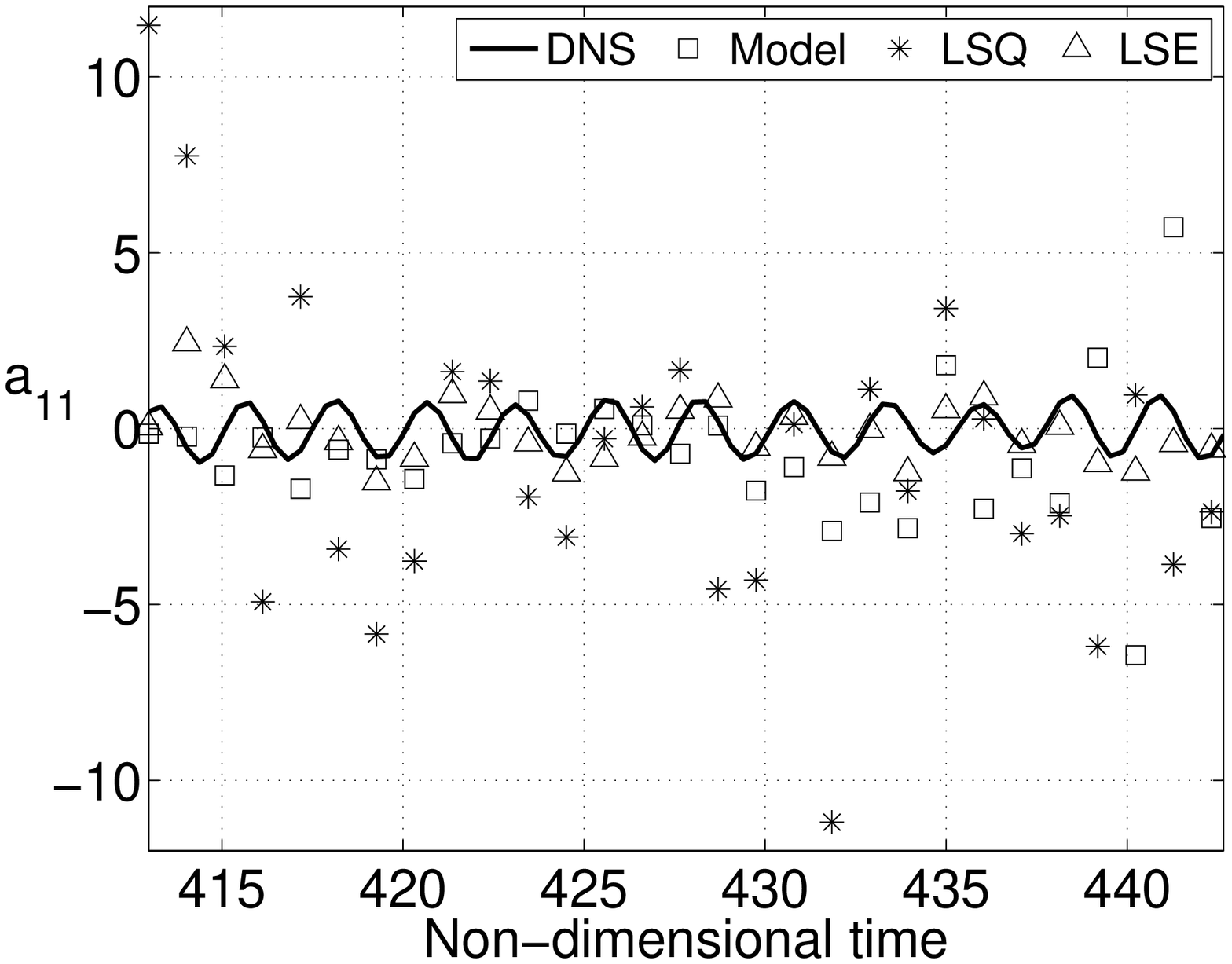}%
  \hspace{1.5cm}
  \includegraphics[width=5.cm,height=4.25cm]{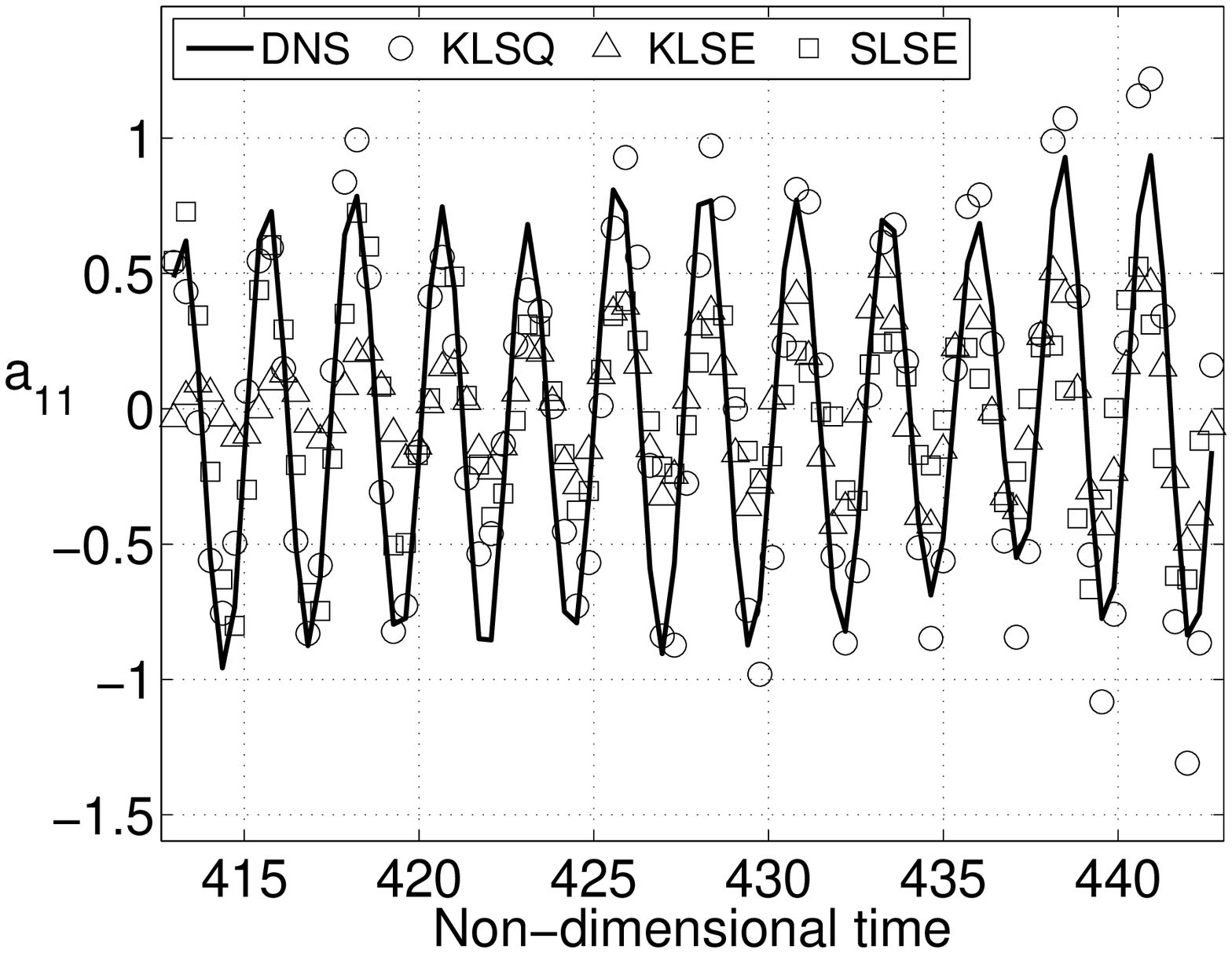}\\%
  \includegraphics[width=5.cm,height=4.25cm]{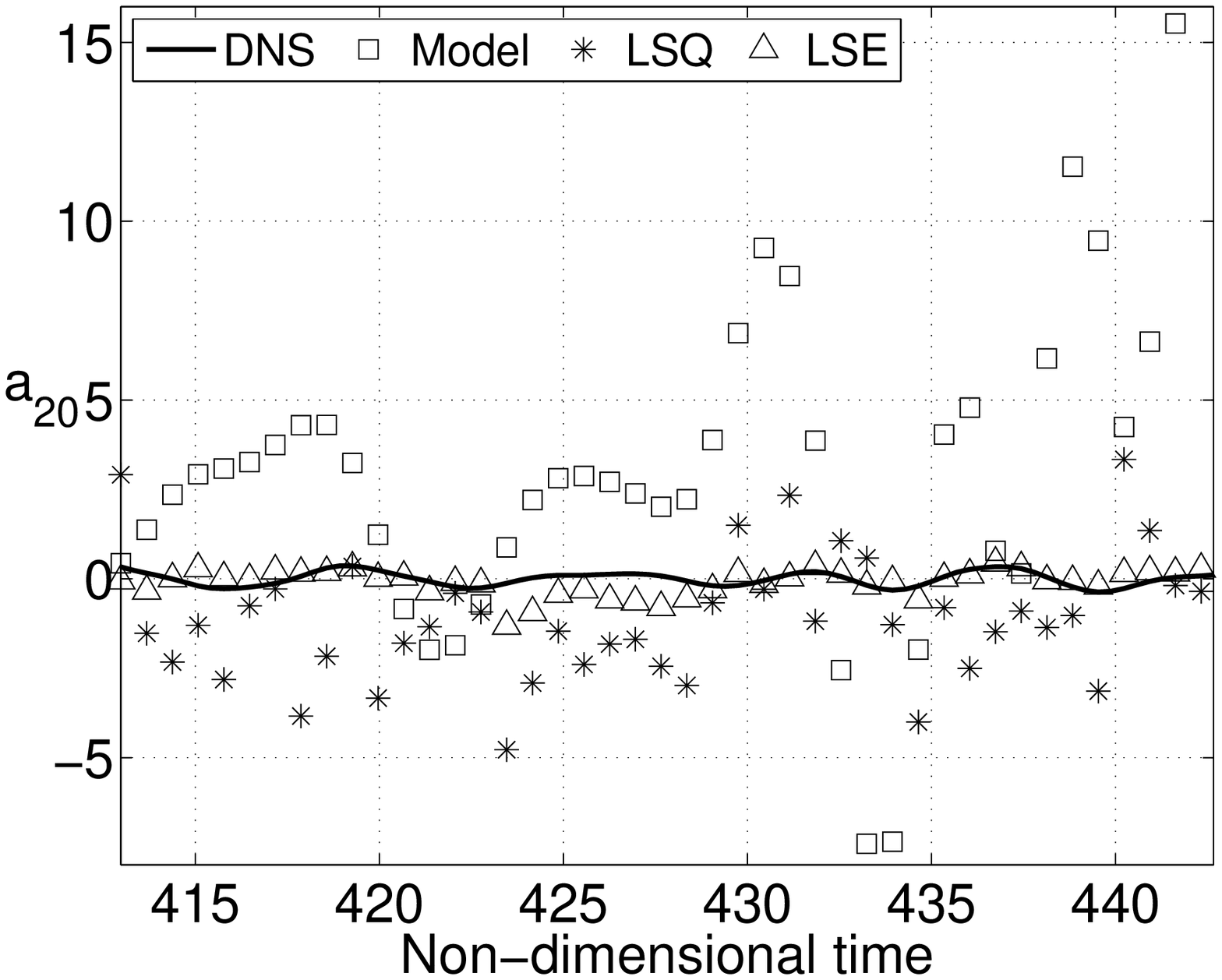}%
  \hspace{1.5cm}
  \includegraphics[width=5.cm,height=4.25cm]{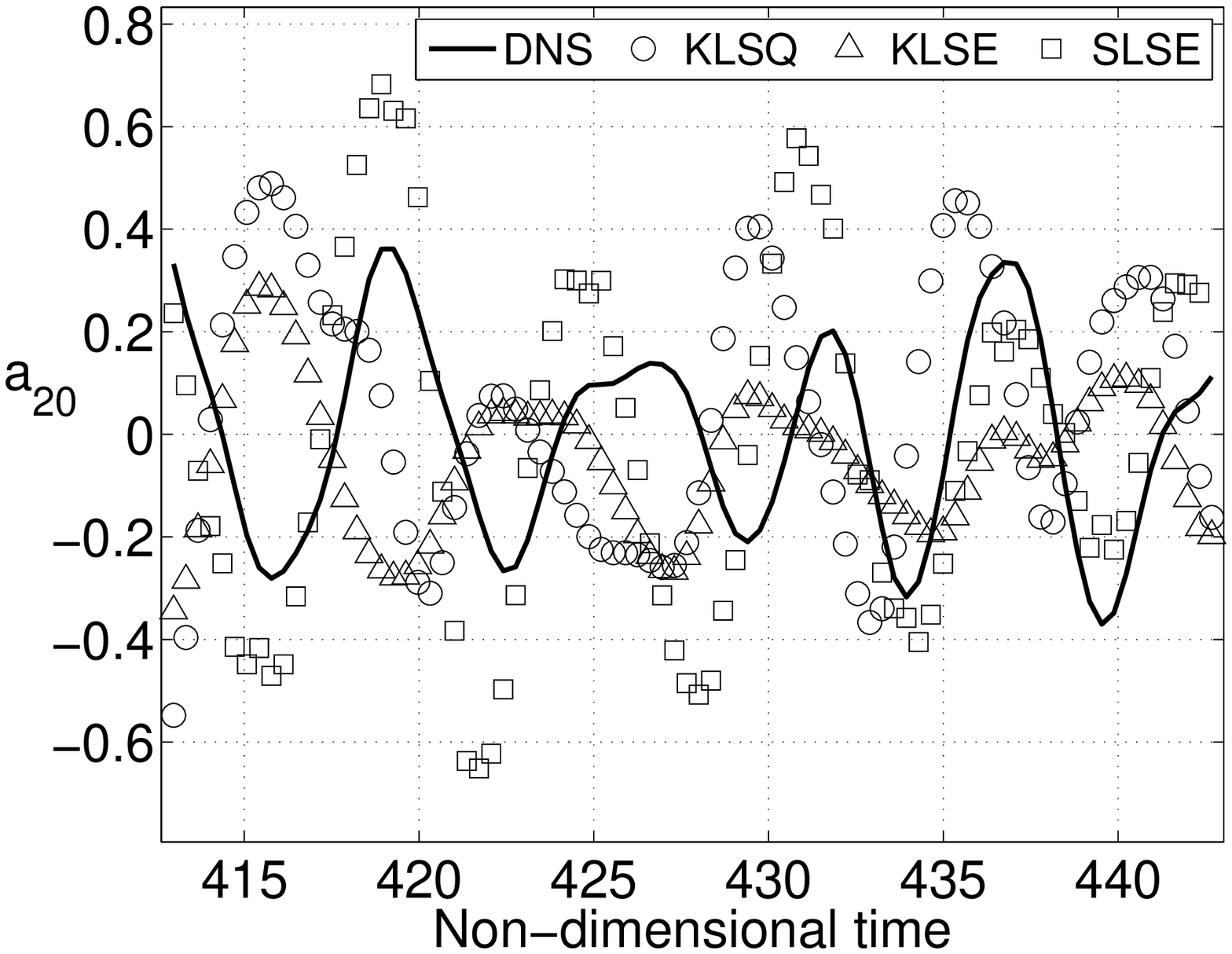}%
  \caption{Estimation of some representative modal coefficients in the
    three-dimensional case, for the configuration (b) (see text) and
    in the time interval close to the calibration one,
    together with the reference ones evaluated
    from the DNS simulation. Note that different axis scales are used
    on the left and on the right plots.}
  \label{fig:coeff3d}
\end{figure}

\begin{table}[h!]
  \centering
  \subfigure{\begin{tabular}{|c|c|c|c|}
    \hline
    SC (a) & $\overline{e(U)}\%$ & $\overline{e(V)}\%$ &
    $\overline{e(W)}\%$ \\%
    \hline
    KLSQ  & 9.85  & 37.00 & 108.12 \\%
    KLSE  & 9.60  & 36.49 & 103.92 \\%
    \hline
    & $\overline{e(U')}\%$ & $\overline{e(V')}\%$ &
    $\overline{e(W')}\%$\\%
    \hline
    KLSQ  & 61.89 & 52.48 & 105.76 \\%
    KLSE  & 9.60  & 36.49 & 103.92 \\%
    \hline
    & $\overline{e(U_f)}\%$ & $\overline{e(V_f)}\%$ &
    $\overline{e(W_f)}\%$ \\%
    \hline
    KLSQ  & 5.08  & 19.40 & 122.54 \\%
    KLSE  & 3.78  & 15.92 &  86.68 \\%
    \hline
  \end{tabular}}
\hspace{5.mm}
  \subfigure{\begin{tabular}{|c|c|c|c|}
    \hline
    SC (b) & $\overline{e(U)}\%$ & $\overline{e(V)}\%$ &
    $\overline{e(W)}\%$ \\%
    \hline
     KLSQ  &     10.00  &     36.65  &    111.52  \\%
     KLSE  &     10.03  &     38.19  &    105.57  \\%
    \hline
    & $\overline{e(U')}\%$ & $\overline{e(V')}\%$ &
    $\overline{e(W')}\%$\\%
    \hline
     KLSQ  &     62.76  &     51.95  &    109.13  \\%
     KLSE  &     62.93  &     54.10  &    103.27  \\%
    \hline
    & $\overline{e(U_f)}\%$ & $\overline{e(V_f)}\%$ &
    $\overline{e(W_f)}\%$ \\%
    \hline
     KLSQ  &     4.18  &     15.87  &     96.63 \\%
     KLSE  &     3.61  &     16.50  &     72.35 \\%
    \hline
  \end{tabular}}
  \subfigure{\begin{tabular}{|c|c|c|c|}
    \hline
    SC (c) & $\overline{e(U)}\%$ & $\overline{e(V)}\%$ &
    $\overline{e(W)}\%$ \\%
    \hline
     KLSQ  &     10.36  &     38.47  &    112.72  \\%
     KLSE  &     10.42  &     38.59  &    108.99  \\%
    \hline
    & $\overline{e(U')}\%$ & $\overline{e(V')}\%$ &
    $\overline{e(W')}\%$\\%
    \hline
     KLSQ  &     65.02  &     54.55  &    110.30  \\%
     KLSE  &     65.37  &     54.70  &    106.68  \\%
    \hline
    & $\overline{e(U_f)}\%$ & $\overline{e(V_f)}\%$ &
    $\overline{e(W_f)}\%$ \\%
    \hline
     KLSQ  &     5.61  &     22.97  &    137.30 \\%
     KLSE  &     4.22  &     18.51  &     97.95 \\%
    \hline
  \end{tabular}}
\hspace{5.mm}
  \subfigure{\begin{tabular}{|c|c|c|c|}
    \hline
    SC (d) & $\overline{e(U)}\%$ & $\overline{e(V)}\%$ &
    $\overline{e(W)}\%$ \\%
    \hline
     KLSQ  &     10.21  &     37.75  &    109.98  \\%
     KLSE  &     10.37  &     38.19  &    108.10  \\%
    \hline
    & $\overline{e(U')}\%$ & $\overline{e(V')}\%$ &
    $\overline{e(W')}\%$\\%
    \hline
     KLSQ  &     64.15  &     53.54  &    107.56  \\%
     KLSE  &     65.08  &     54.12  &    105.82  \\%
    \hline
    & $\overline{e(U_f)}\%$ & $\overline{e(V_f)}\%$ &
    $\overline{e(W_f)}\%$ \\%
    \hline
     KLSQ  &     5.02  &     20.16  &    121.39 \\%
     KLSE  &     4.03  &     16.41  &     91.39 \\%
    \hline
  \end{tabular}}
  \subfigure{\begin{tabular}{|c|c|c|c|}
    \hline
    SC (e) & $\overline{e(U)}\%$ & $\overline{e(V)}\%$ &
    $\overline{e(W)}\%$ \\%
    \hline
     KLSQ  &     11.30  &     40.99  &    119.15  \\%
     KLSE  &     10.24  &     38.35  &    107.59  \\%
    \hline
    & $\overline{e(U')}\%$ & $\overline{e(V')}\%$ &
    $\overline{e(W')}\%$\\%
    \hline
     KLSQ  &     71.01  &     58.17  &    116.48  \\%
     KLSE  &     64.23  &     54.31  &    105.33  \\%
    \hline
    & $\overline{e(U_f)}\%$ & $\overline{e(V_f)}\%$ &
    $\overline{e(W_f)}\%$ \\%
    \hline
     KLSQ  &     6.78  &     25.30  &    166.46 \\%
     KLSE  &     3.70  &     16.61  &     87.30 \\%
    \hline
  \end{tabular}}
  \caption{Relative percent errors, in $L^2$ norm, in the
    reconstruction of the velocity components ($\overline{e(U)}$,
    $\overline{e(V)}$, $\overline{e(W)}$), of their fluctuating part
    ($\overline{e(U')}$, $\overline{e(V')}$, $\overline{e(W')}$), and
    of the part of the velocity field projected on the retained POD
    modes ($\overline{e(U_f)}$, $\overline{e(V_f)}$,
    $\overline{e(W_f)}$). The errors are averaged in time in the
    estimation interval, which starts just after the calibration
    interval. SC: Sensor Configuration.}  
  \label{tab:ric3d_close}
\end{table}
The parameter $C_R$ of equation~\ref{klsq} was selected by
experimenting different values. For example, in figure~\ref{fig:CR} we show the
$L^2$ relative error in the
reconstruction of the U component
projected on the retained POD modes, as a function of $C_R$. The
results are relative to the configuration (b) 
in the time interval close to the
calibration one, using K-LSE.
For all of the results shown in the
following we
took $C_R=0.1$ for K-LSE. Note that for $C_R \ge 10^2$ the results are
basically those of a simple LSE. A similar analysis was performed using
K-LSQ and the optimal value that was selected is $C_R=10$.
Results relative to the configuration (b) 
in the time interval close to the
calibration one 
are reported in figure~\ref{fig:coeff3d}, 
where some representative modal coefficients
predicted by the 
calibrated POD model, LSQ, LSE,  K-LSQ, K-LSE  and SLSE are
plotted, together with the projection of the DNS velocity fields on
the corresponding POD mode. 

\begin{figure}[h!]
   \centering
   \begin{tabular}{ccc}
     \subfigure{\includegraphics[width=4.5cm,height=3.75cm]{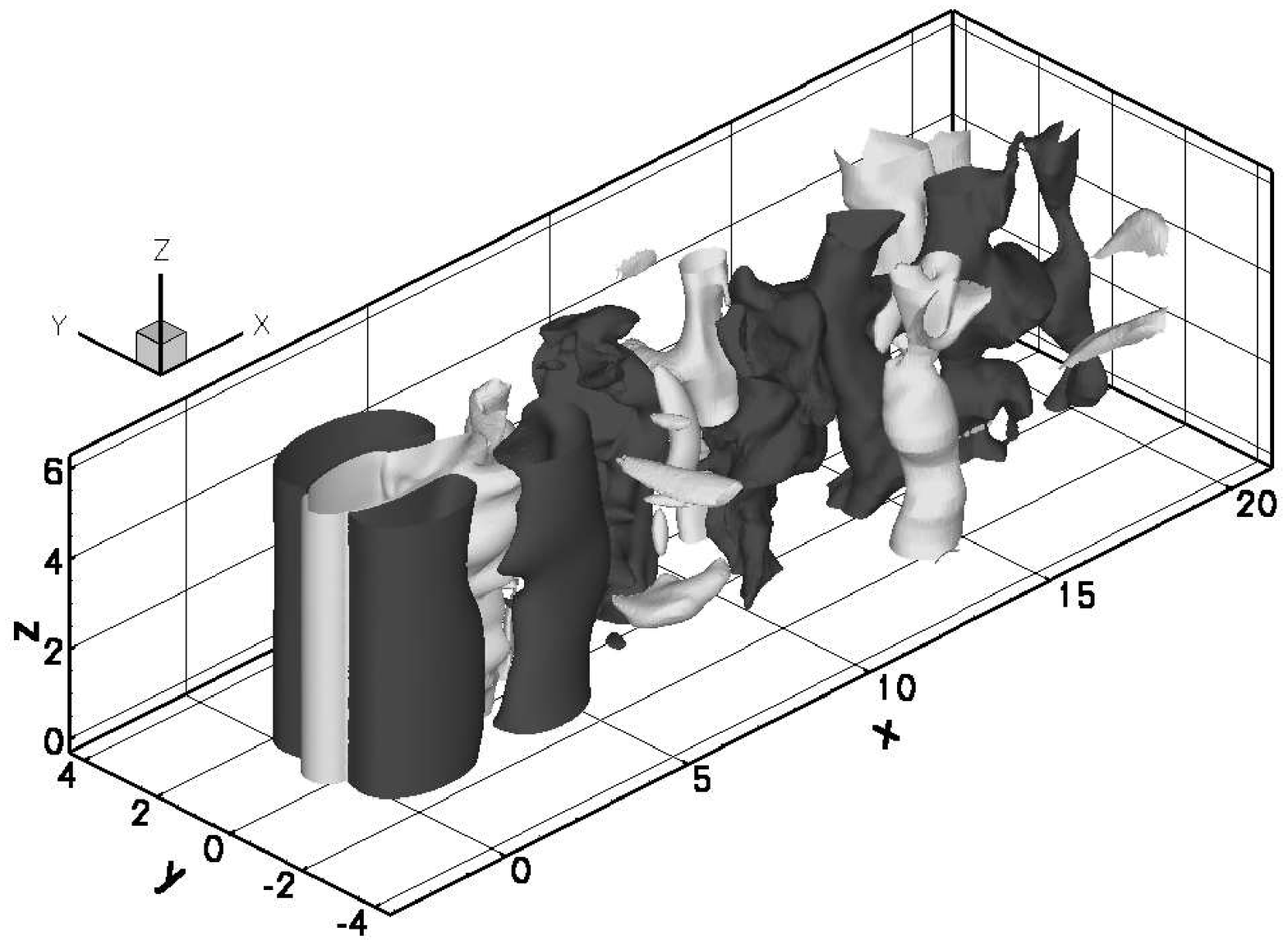}}%
     & \hspace{-3.5mm}
     \subfigure{\includegraphics[width=4.5cm,height=3.75cm]{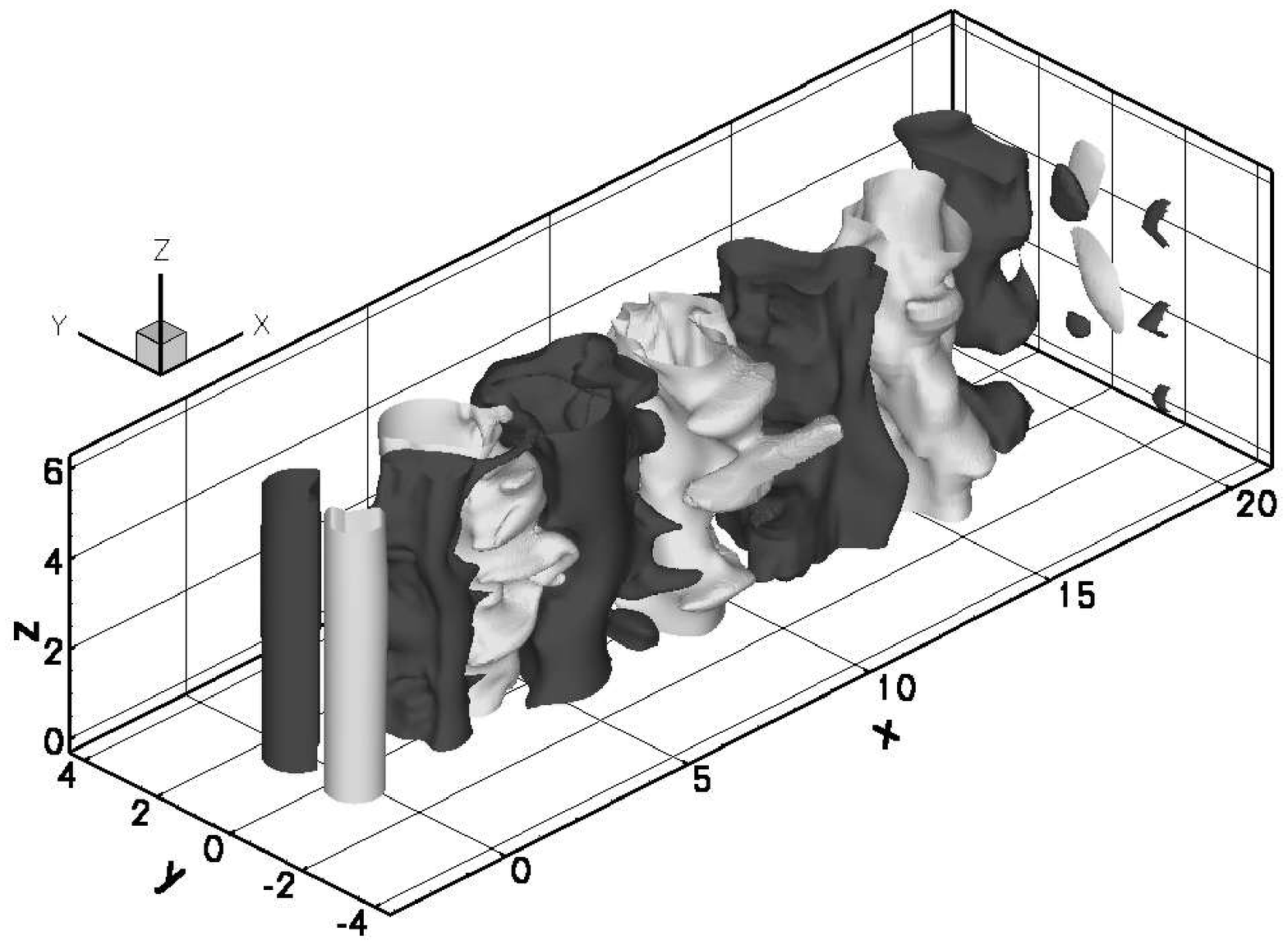}}%
     & \hspace{-3.5mm}
     \subfigure{\includegraphics[width=4.5cm,height=3.75cm]{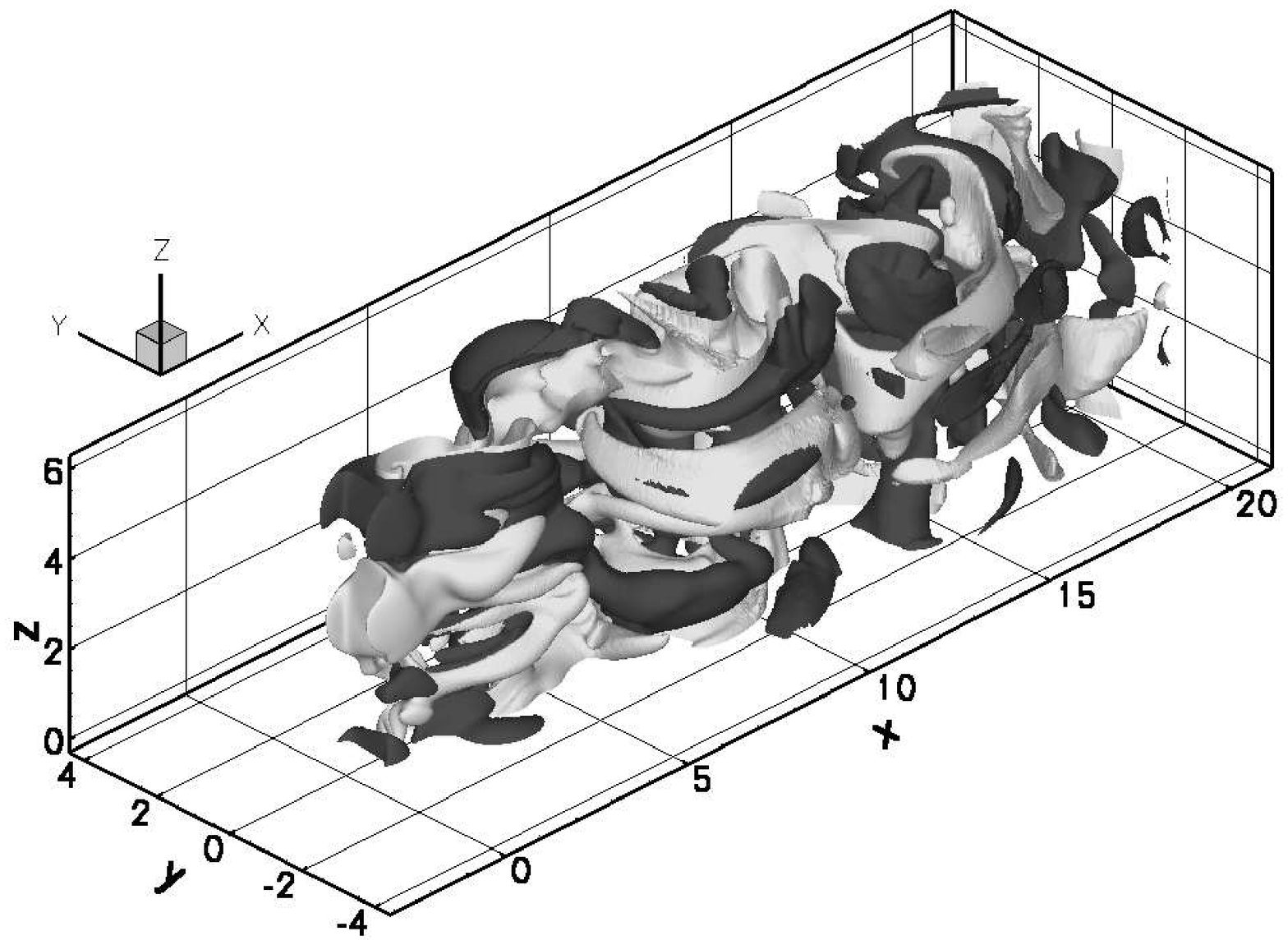}}\\
     & (a) & \\
     \subfigure{\includegraphics[width=4.5cm,height=3.75cm]{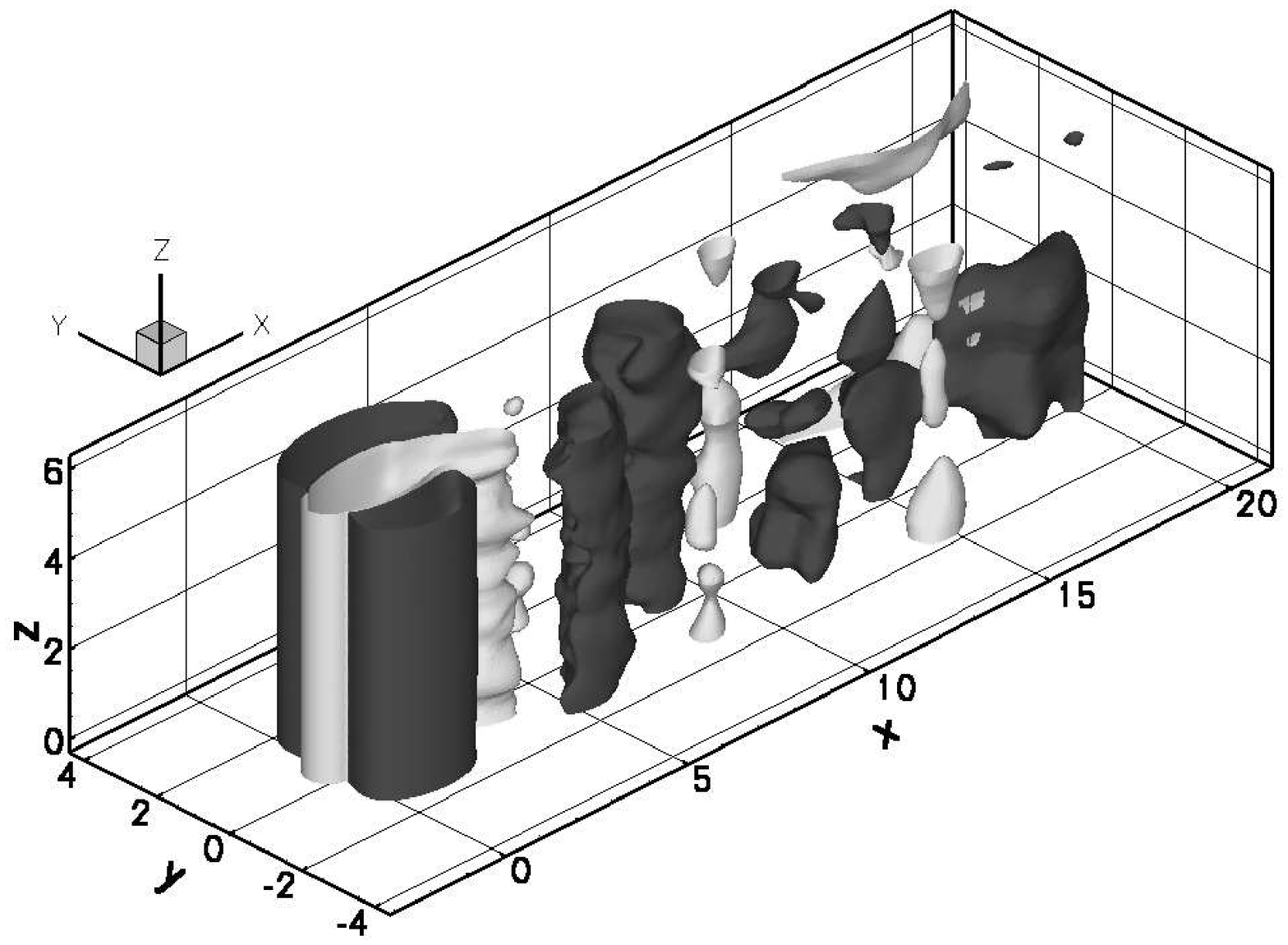}}%
     & \hspace{-3.5mm}
     \subfigure{\includegraphics[width=4.5cm,height=3.75cm]{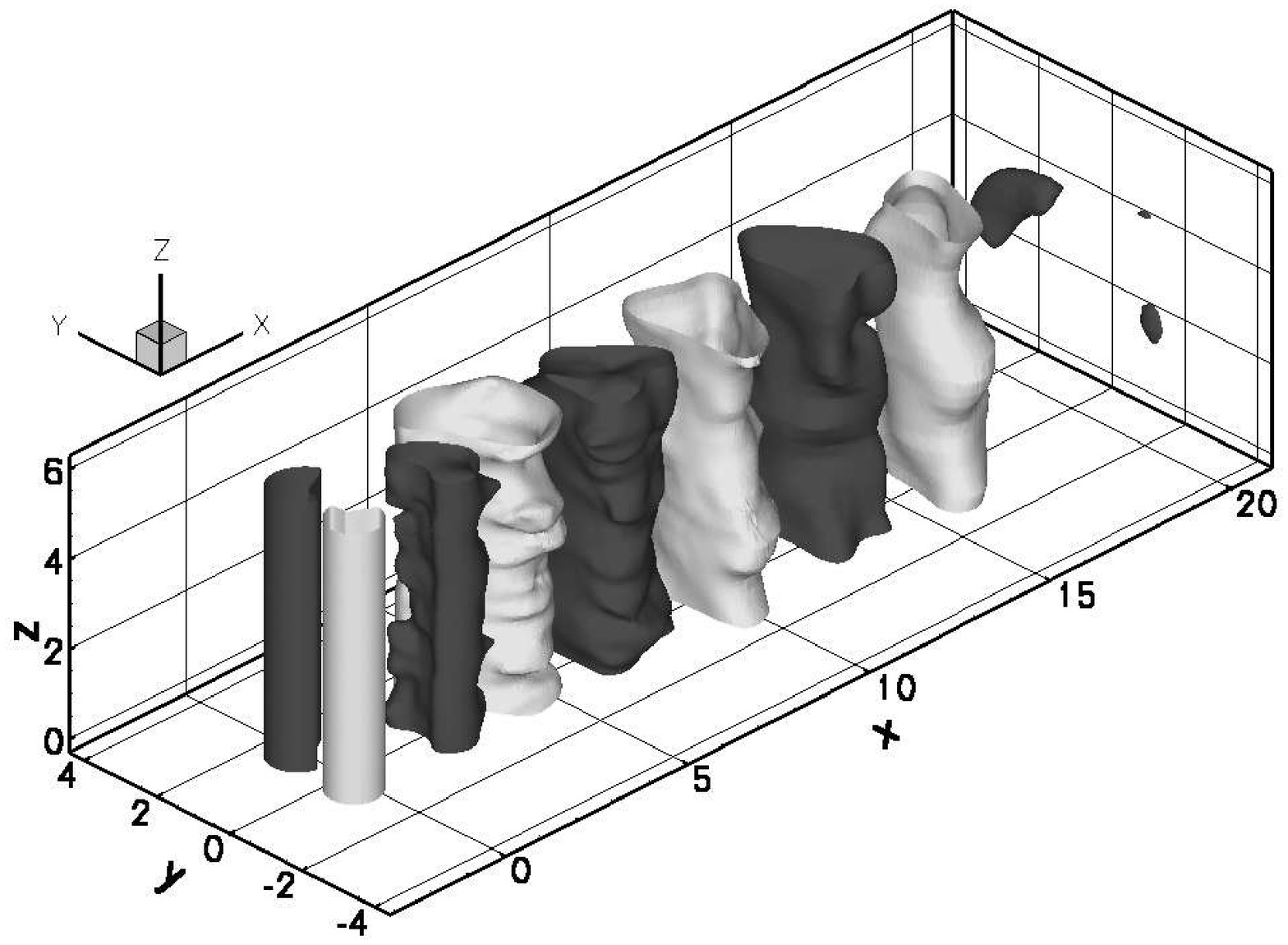}}%
    & \hspace{-3.5mm}
    \subfigure{\includegraphics[width=4.5cm,height=3.75cm]{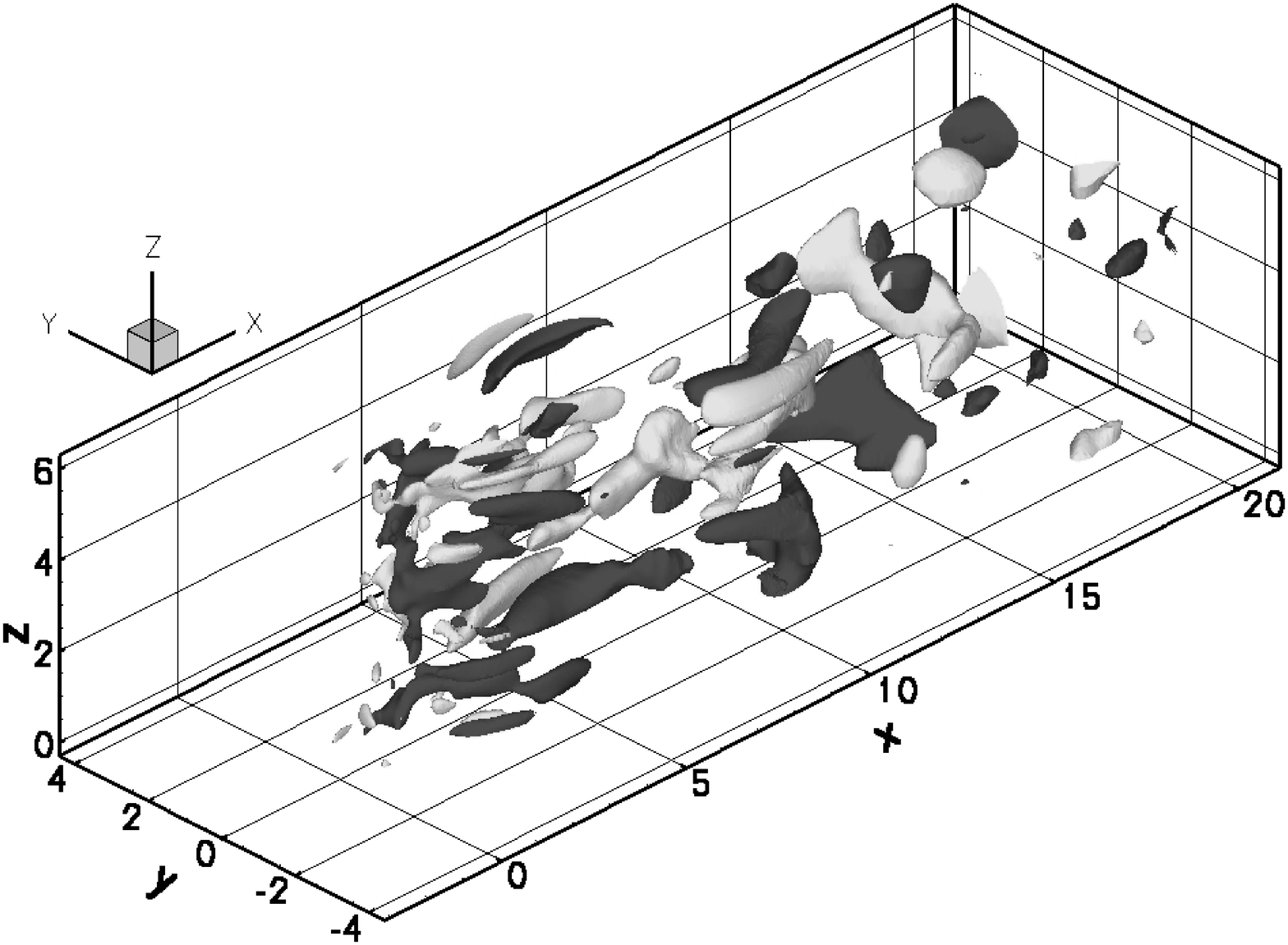}}\\
    & (b) & \\
    \subfigure{\includegraphics[width=4.5cm,height=3.75cm]{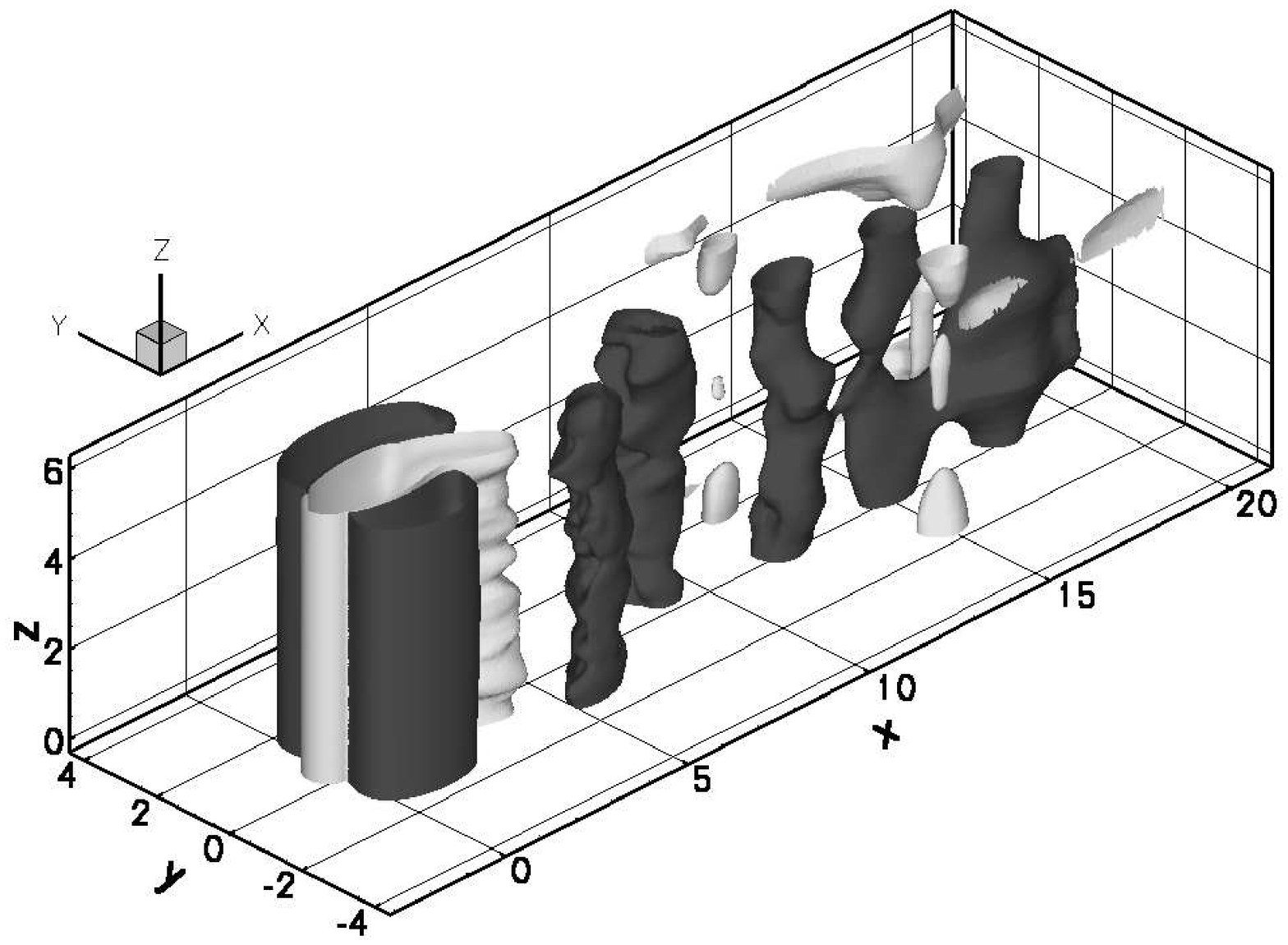}}%
    & \hspace{-3.5mm}
    \subfigure{\includegraphics[width=4.5cm,height=3.75cm]{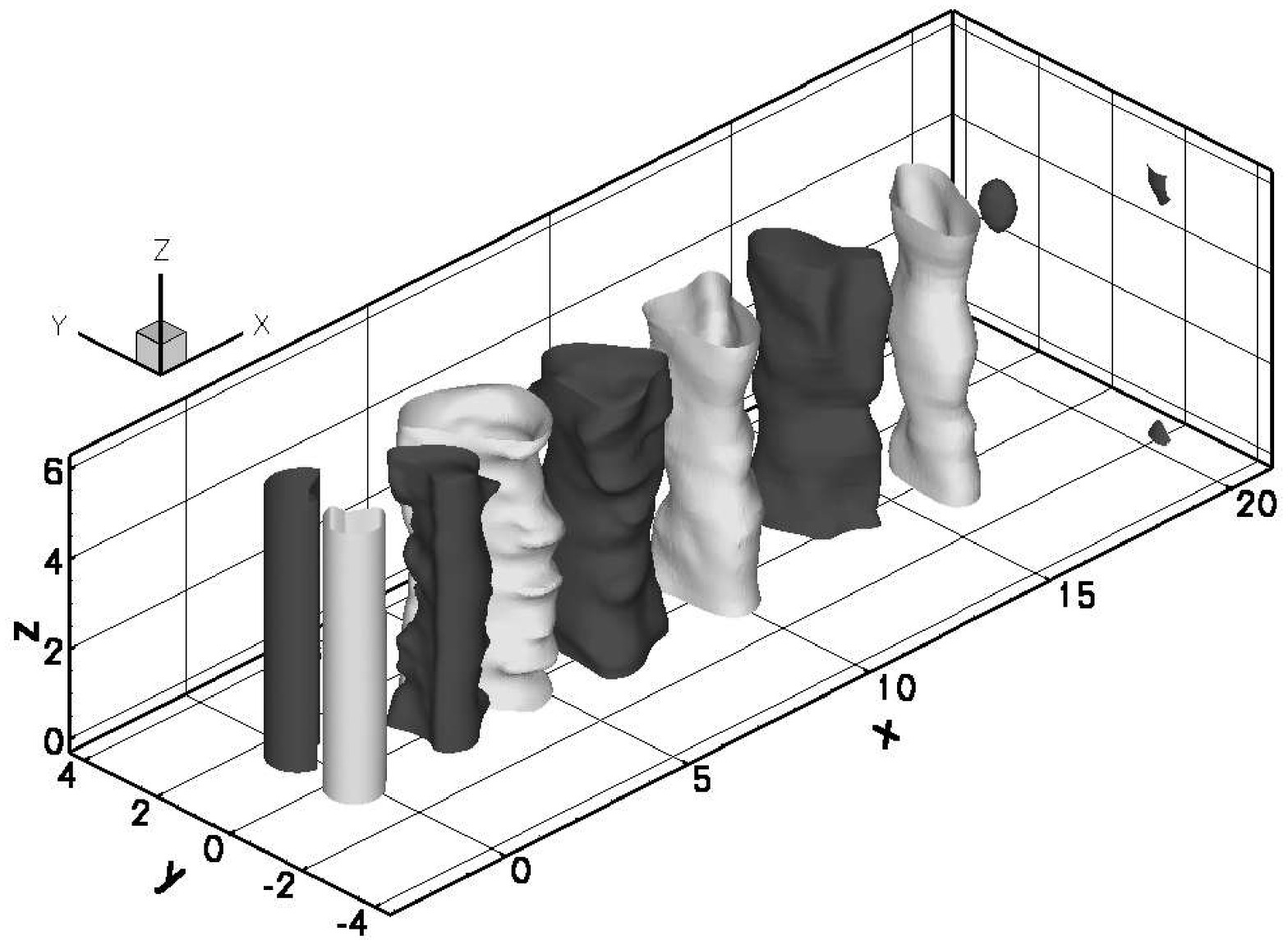}}%
    & \hspace{-3.5mm}
    \subfigure{\includegraphics[width=4.5cm,height=3.75cm]{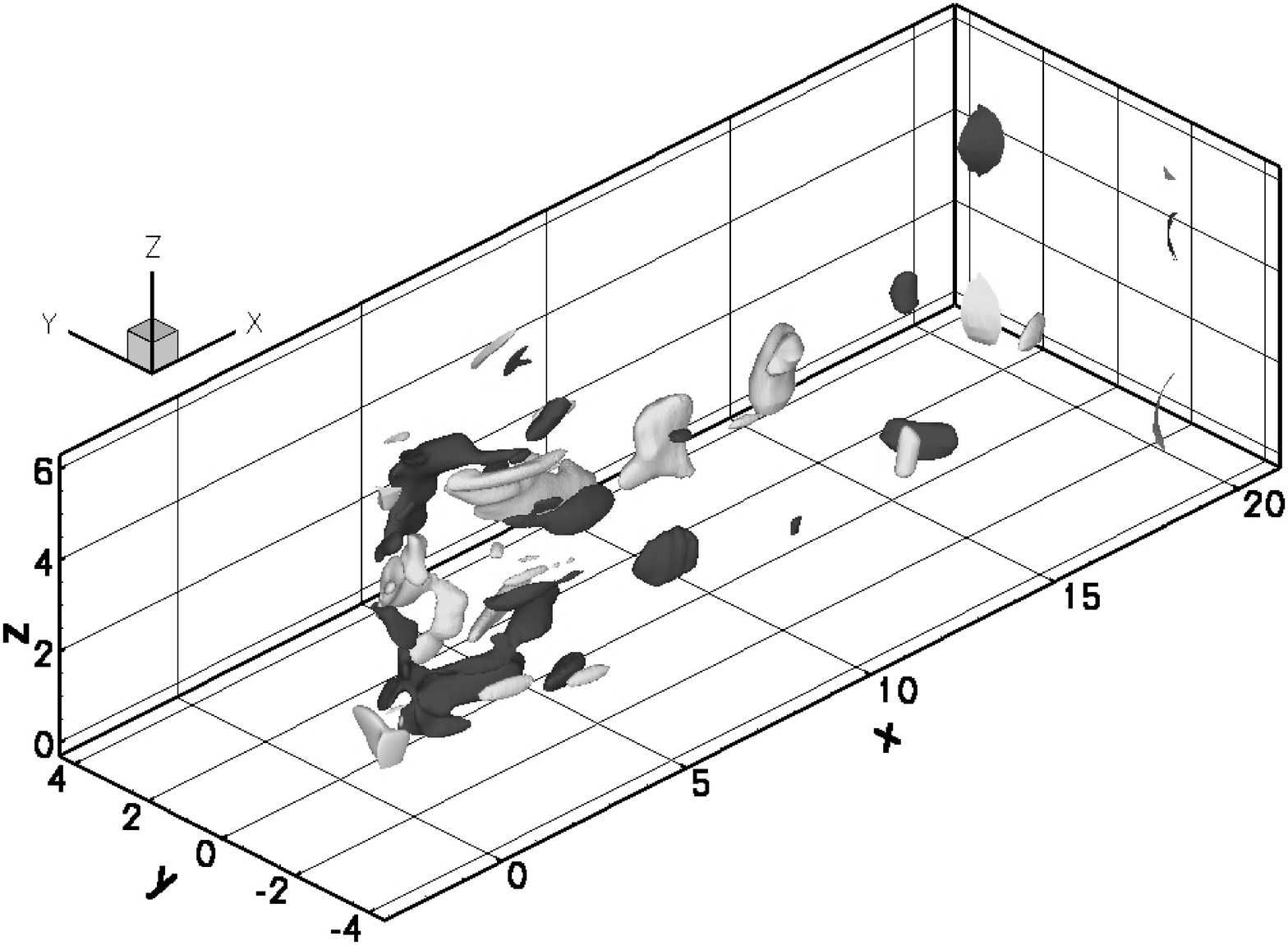}}\\
    & (c) & \\
  \end{tabular}
  \caption{Isosurfaces of the velocity components $u$ (left, grey =
    0.5, dark grey = 
    1.0), $v$ (center, grey = -0.25, dark grey = 0.25) and $w$ 
    (right, grey = -0.075, dark grey = 0.075) of a snapshot outside
    the database: (a) actual snapshot, (b) snapshot projected on the
    retained POD modes,(c) reconstructed snapshot using the
    K-LSE technique with the sensor configuration (b).}
  \label{fig:snapric}
\end{figure}

Results of configuration (b) has been shown because the placement of
the sensors is appropriate for the LSE method, as already discussed,
and this makes the comparison with the proposed approaches
more comprehensive.
\begin{table}[h!]
  \centering
  \subfigure{\begin{tabular}{|c|c|c|c|}
    \hline
    SC (a) & $\overline{e(U)}\%$ & $\overline{e(V)}\%$ &
    $\overline{e(W)}\%$ \\%
    \hline
    KLSQ  &      8.79  &     32.56  &    114.51 \\%
    KLSE  &      8.96  &     33.54  &    110.95  \\%
    \hline
    & $\overline{e(U')}\%$ & $\overline{e(V')}\%$ &
    $\overline{e(W')}\%$\\%
    \hline
    KLSQ  &      55.91  &     46.13  &    111.51  \\%
    KLSE  &      56.97  &     47.53  &    108.00  \\%
    \hline
    & $\overline{e(U_f)}\%$ & $\overline{e(V_f)}\%$ &
    $\overline{e(W_f)}\%$ \\%
    \hline
    KLSQ  &      4.13  &     15.74  &    101.87 \\%
    KLSE  &      3.63  &     14.90  &     84.19 \\%
    \hline
  \end{tabular}}
\hspace{5.mm}
  \subfigure{\begin{tabular}{|c|c|c|c|}
    \hline
    SC (b) & $\overline{e(U)}\%$ & $\overline{e(V)}\%$ &
    $\overline{e(W)}\%$ \\%
    \hline
    KLSQ  &      9.50  &     34.47  &    121.52\\%
    KLSE  &      8.83  &     33.11  &    111.35\\%
    \hline
    & $\overline{e(U')}\%$ & $\overline{e(V')}\%$ &
    $\overline{e(W')}\%$\\%
    \hline
    KLSQ  &      60.50  &     48.87  &    118.35  \\%
    KLSE  &      56.17  &     46.91  &    108.42  \\%
    \hline
    & $\overline{e(U_f)}\%$ & $\overline{e(V_f)}\%$ &
    $\overline{e(W_f)}\%$ \\%
    \hline
    KLSQ  &      5.15  &     19.82  &    118.19 \\%
    KLSE  &      3.34  &     14.68  &     80.84 \\%
    \hline
  \end{tabular}}
  \subfigure{\begin{tabular}{|c|c|c|c|}
    \hline
    SC (c) & $\overline{e(U)}\%$ & $\overline{e(V)}\%$ &
    $\overline{e(W)}\%$ \\%
    \hline
    KLSQ  &      9.47  &     33.88  &    120.45\\%
    KLSE  &      9.17  &     33.19  &    110.95\\%
    \hline
    & $\overline{e(U')}\%$ & $\overline{e(V')}\%$ &
    $\overline{e(W')}\%$\\%
    \hline
    KLSQ  &      60.08  &     47.97  &    117.20  \\%
    KLSE  &      58.36  &     47.03  &    108.09  \\%
    \hline
    & $\overline{e(U_f)}\%$ & $\overline{e(V_f)}\%$ &
    $\overline{e(W_f)}\%$ \\%
    \hline
    KLSQ  &      5.28  &     19.02  &    133.11 \\%
    KLSE  &      3.07  &     13.56  &     74.63 \\%
    \hline
  \end{tabular}}
\hspace{5.mm}
  \subfigure{\begin{tabular}{|c|c|c|c|}
    \hline
    SC (d) & $\overline{e(U)}\%$ & $\overline{e(V)}\%$ &
    $\overline{e(W)}\%$ \\%
    \hline
    KLSQ  &      9.42  &     35.45  &    121.86  \\%
    KLSE  &      9.44  &     34.40  &    110.40  \\%
    \hline
    & $\overline{e(U')}\%$ & $\overline{e(V')}\%$ &
    $\overline{e(W')}\%$\\%
    \hline
    KLSQ  &      59.93  &    50.26  &    118.64  \\%
    KLSE  &      59.93  &    48.71  &    107.54  \\%
    \hline
    & $\overline{e(U_f)}\%$ & $\overline{e(V_f)}\%$ &
    $\overline{e(W_f)}\%$ \\%
    \hline
    KLSQ  &      5.11  &     20.46  &    133.56 \\%
    KLSE  &      3.62  &     15.74  &     76.05 \\%
    \hline
  \end{tabular}}
  \subfigure{\begin{tabular}{|c|c|c|c|}
    \hline
    SC (e) & $\overline{e(U)}\%$ & $\overline{e(V)}\%$ &
    $\overline{e(W)}\%$ \\%
    \hline
    KLSQ  &     10.60  &     38.50  &    136.53  \\%
    KLSE  &      9.22  &     33.49  &    110.24  \\%
    \hline
    & $\overline{e(U')}\%$ & $\overline{e(V')}\%$ &
    $\overline{e(W')}\%$\\%
    \hline
    KLSQ  &     67.42  &     54.59  &    132.86  \\%
    KLSE  &     58.65  &     47.45  &    107.38  \\%
    \hline
    & $\overline{e(U_f)}\%$ & $\overline{e(V_f)}\%$ &
    $\overline{e(W_f)}\%$ \\%
    \hline
    KLSQ  &     7.29  &     26.00  &    183.25 \\%
    KLSE  &     3.27  &     14.10  &     77.47 \\%
    \hline
  \end{tabular}}
  \caption{Relative percent errors, in $L^2$ norm, in the
    reconstruction of the velocity components ($\overline{e(U)}$,
    $\overline{e(V)}$, $\overline{e(W)}$), of their fluctuating part
    ($\overline{e(U')}$, $\overline{e(V')}$, $\overline{e(W')}$), and
    of the part of the velocity field projected on the retained POD
    modes ($\overline{e(U_f)}$, $\overline{e(V_f)}$,
    $\overline{e(W_f)}$). The errors are averaged in time in the
    estimation interval, which is about $52$ time units far from the
    calibration interval.}
  \label{tab:ric3d_far}
\end{table}
The results obtained with the other sensor configurations and on
the time interval far from the calibration one are
qualitatively analogous and quantitatively similar to those reported
in figure~\ref{fig:coeff3d}, except for the LSE and the LSQ methods which
are more sensitive to sensor placement. 
In figure~\ref{fig:coeff3d} it is seen that
LSE, LSQ, and also the calibrated POD model
provide reasonable predictions only for the first two modal
coefficients, that are associated with the vortex-shedding
dynamics. The predictions of the remaining
modes are completely unreliable.
When dynamic estimation is applied, or when the SLSE
approach is used, predictions are definitely
improved. 
In particular, this is true for modes like $a_1$ or $a_{11}$ that are
related to the vortex-shedding, 
i.e., almost periodic with a period that is the same or a multiple of
the vortex-shedding 
period. 
The prediction of the remaining modes is definitely less accurate 
(see $a_3$ and $a_{20}$), especially when
very low frequencies are dominant, as in the case of $a_3$. 
However, the overall accuracy is significantly improved in
comparison with the LSQ and LSE approaches alone. 
The K-LSE approach is systematically more accurate than 
SLSE. This is quantitatively confirmed by the relative error
between the DNS velocity components and those reconstructed by
the dynamic approaches and the SLSE one,
reported in table~\ref{tab:ric3d_close} for all the sensor
configurations. 
However, the table shows that
there are no significant differences in accuracy
between the K-LSQ, K-LSE and SLSE methods, and the reconstruction errors
can be considered satisfactory if we
consider the complexity of the flow with respect to the 2D
case. 
Indeed, the streamwise and vertical velocity components are
reconstructed with errors of the order of $9\%$ and $33\%$,
respectively. Note that errors are computed over the whole
computational domain, which extends for $20~L$ behind the
cylinder, and, as it is qualitatively shown in figure \ref{fig:snapric}, the
reconstruction errors are small in the near-wake and they
progressively increase moving away from the cylinder in the downstream
direction. 
Concerning the spanwise component of the velocity, errors are large
since the retained POD modes themselves poorly represent this
component of velocity as it is not energetically significant, in
average, with respect to the remaining ones. This aspect might be
improved working on the construction of the POD basis choosing, for
instance, a different norm which weights more the spanwise component
of the velocity or which corresponds to a quantity different from kinetic
energy. 
The K-LSQ, K-LSE and SLSE methods are also similar in the
sensitivity of the predictions to sensors type and placement, which is
generally low. Nevertheless, the predictions given by the
K-LSE method are systematically  the most
insensitive to sensors placement. 
Table~\ref{tab:ric3d_close} also shows the
relative reconstruction errors on the fluctuating part of the velocity
components. 
The errors in the reconstruction of the fluctuating part of both the
streamwise and the lateral velocity components are comparable.
To help the interpretation of the errors in the reconstruction of the
fluctuating flow field, which might seem large at first sight,
it is important to point out that, outside the
calibration interval, the amount of fluctuating energy that can be
recovered by 20 modes is slightly more than $60\%$
(\cite{Buffoni2006}). Therefore, the error in the velocity
components is due in part to a pure approximation error. In other
words, the
accuracy of the best possible reconstruction is limited from above by
the capability of the POD modes to actually represent the flow outside
the time interval where the snapshots were taken, which 
however, increases using a larger snapshots database, as
shown in  \cite{Buffoni2006}.

In order to separate the accuracy of the estimation methods from the
representativeness of the POD basis,
we  evaluated 
the errors between the estimated field and  the 
projections of the velocity fields on the retained POD modes. Such
errors are 
reported in tables~\ref{tab:ric3d_close} and \ref{tab:ric3d_far}.
These errors are systematically and significantly lower
than those observed in the reconstruction of the whole flow fields.

Comparing the results reported in table~\ref{tab:ric3d_close} with
those of table~\ref{tab:ric3d_far}, where the same quantities are
shown for the time interval far from the calibration one, it can be
seen that
the accuracy of the predictions is independent of the distance from the
calibration interval, this being a very positive feature of all the
proposed approaches. 

In figure~\ref{fig:snapric} the velocity components obtained
by DNS at $t=426.6$ (a snapshot outside the database used
for the derivation and calibration of the POD model) are plotted
together with their projection  in the space of the retained POD basis, which
represents the best approximation of the flow which can be estimated
with the retained POD modes,
and with the prediction given by the K-LSE method. 
It can be seen that the main  structures 
characterizing the streamwise and lateral velocity fields are well
reconstructed. As for
the spanwise velocity component, the reconstruction accuracy is not
satisfactory, but this is due to the fact that it is one order of
magnitude lower than the other components, as already
discussed. 

\section{Conclusions}

We devised a method to construct a non-linear observer for unsteady 
flows. This method is based on the coupling of a non-linear
low-dimensional model of the flow with a linear technique that 
estimates the
coefficients of the flow representation in terms of POD modes. 
The underlying idea is that the estimated flow should
approximately satisfy the POD model. 
The coupling leads to
a nonlinear minimization problem solved by a pseudo-spectral approach
and a Newton method.   

The non-linear observer was applied to the laminar flow around a
confined square cylinder at two different Reynolds numbers; at the
first the flow is two-dimensional, while in the second case
three-dimensional phenomena occur in the wake. 
In the two-dimensional case, since
the flow patterns are rather simple, the results show that the
proposed procedure is able to 
give a significantly more accurate estimate of the POD coefficients
even with a limited number of sensors, than those obtained with
the LSQ and the LSE approaches.
The QSE method improves the predictions of the LSE, even if the
overall accuracy remains definitely lower than those given by
the proposed approaches. 
In the three-dimensional case, the flow 
dynamics is more complex, and not only LSE and  LSQ, but also the
calibrated POD dynamical system provide unsatisfactory coefficient
estimations when used outside the calibration interval. 
Conversely, the proposed procedure, combined with either LSQ
or LSE, 
gives accurate predictions of the coefficients of those POD
modes that are related to vortex shedding. For the remaining modes,
the accuracy is lower. 
Nevertheless, 
the instantaneous velocity fields are reconstructed with
satisfactory accuracy, both close and far from the calibration
interval of the POD model. 
Moreover, K-LSE and K-LSQ methods are weakly sensitive to
sensor type and placement.
The results obtained with the proposed approaches are comparable to
those obtained by the SLSE approach, which also uses the temporal
history of the flow measurements, but in the Fourier space. However,
this latter techniques has computational complexity which is
significantly larger than that of LSE or LSQ, and it is comparable to
that of the proposed approaches.

\vskip.2cm

This work was funded in part by the HPC-EUROPA project
(RII3-CT-2003-506079). 
IDRIS (Orsay, France) and M3PEC
(Universit\'e Bordeaux 1, France) provided the
computational resources.

\bibliography{rr_bcils}
\bibliographystyle{plain}
\end{document}